\documentclass[aos, authoryear]{imsart}
\RequirePackage[OT1]{fontenc}
\RequirePackage{amsthm,amsmath,natbib}
\RequirePackage{hypernat}
\startlocaldefs
\numberwithin{equation}{section}
\theoremstyle{plain}

\endlocaldefs
\usepackage{psfrag, epsfig, graphicx}
\usepackage{amsfonts,amsmath,latexsym,amssymb}
\usepackage[active]{srcltx}
\usepackage{fullpage}

\newtheorem{lemma}{Lemma}
\newtheorem{theorem}{Theorem}
\newtheorem{proposition}{Proposition}

\newtheorem{assumption}{Assumption}

\newtheorem*{ass4.1}{Assumption $4^\prime$}
\newtheorem*{ass4.2}{Assumption $4^{\prime\prime}$}

\newtheorem{remark}{Remark}
\newtheorem{definition}{Definition}
\newtheorem{corollary}{Corollary}

\renewcommand{\kappa}{\varkappa}

\newcommand{\e}{\varepsilon}

\newcommand{\ry}{\mathrm{y}}

\newcommand{\rd}{\mathrm{d}}

\newcommand{\rk}{\mathrm{k}}
\newcommand{\rh}{\mathrm{h}}

\newcommand{\cB}{{\cal B}}
\newcommand{\cC}{{\cal C}}

\newcommand{\cE}{{\cal E}}

\newcommand{\cH}{{\cal H}}
\newcommand{\cI}{{\cal I}}

\newcommand{\cL}{\mathcal{L}}

\newcommand{\cN}{{\cal N}}

\newcommand{\cP}{{\cal P}}
\newcommand{\cQ}{{\cal Q}}
\newcommand{\cR}{{\cal R}}

\newcommand{\cU}{{\cal U}}

\newcommand{\cX}{{\cal X}}
\newcommand{\cY}{{\cal Y}}

\newcommand{\mfi}{\mathbf{i}}
\newcommand{\mfI}{\mathbf{I}}
\newcommand{\mfj}{\mathbf{j}}
\newcommand{\mfJ}{\mathbf{J}}
\newcommand{\bB}{\mathbb B}

\newcommand{\bE}{\mathbb E}

\newcommand{\bL}{{\mathbb L}}

\newcommand{\bN}{{\mathbb N}}

\newcommand{\bP}{{\mathbb P}}

\newcommand{\bR}{{\mathbb R}}

\newcommand{\bX}{{\mathbb X}}

\newcommand{\bZ}{{\mathbb Z}}
\newcommand{\mA}{\mathfrak{A}}

\newcommand{\mF}{\mathfrak{F}}
\newcommand{\mH}{\mathfrak{H}}

\newcommand{\mL}{\mathfrak{L}}

\newcommand{\mP}{\mathfrak{P}}

\newcommand{\mh}{\mathfrak{h}}

\newcommand{\ma}{\mathfrak{a}}

\newcommand{\mb}{\mathfrak{b}}

\newcommand{\epr}{\hfill\hbox{\hskip 4pt
                \vrule width 5pt height 6pt depth 1.5pt}\vspace{0.5cm}\par}

\newcommand{\blh}{\boldsymbol{h}}

\newcommand{\blalpha}{\boldsymbol{\alpha}}
\newcommand{\blb}{\boldsymbol{\beta}}

\newcommand{\blI}{\boldsymbol{I}}
\newcommand{\blP}{\boldsymbol{\mathcal{P}}}
\begin{document}

\begin{frontmatter}
\title{Multivariate density estimation under sup-norm loss: oracle approach, adaptation and independence structure}
\runtitle{Density estimation}
\begin{aug}
\author{\fnms{Oleg} \snm{Lepski}
\ead[label=e2]{lepski@cmi.univ-mrs.fr}}
\runauthor{A. Goldenshluger and O. Lepski}

\affiliation{ Universit\'e Aix--Marseille }


\address{Laboratoire d'Analyse, Topologie, Probabilit\'es\\
 Universit\'e Aix-Marseille  \\
 39, rue F. Joliot-Curie \\
13453 Marseille, France\\
\printead{e2}\\ }
\end{aug}
\begin{abstract}
The paper deals with the density estimation on $\bR^d$ under
sup-norm loss. We provide with fully data-driven estimation procedure
and establish for it  so called sup-norm oracle inequality. The proposed estimator
allows to
take into account not only approximation properties of the underlying density but
 eventual independence structure as well. Our results  contain,
as a particular case, the complete solution of the bandwidth selection
problem in multivariate density model.
Usefulness of the developed approach  is illustrated
by application
to adaptive estimation over anisotropic Nikolskii classes.

\end{abstract}

\begin{abstract}

\end{abstract}


\end{frontmatter}

\def\huh{\hbox{\vrule width 2pt height 8pt depth 2pt}}
\def\blacksquare{\hbox{\vrule width 4pt height 4pt depth 0pt}}
\def\square{\hbox{\vrule\vbox{\hrule\phantom{o}\hrule}\vrule}}
\def\inter{\mathop{{\rm int}}}
%
%
\section{Introduction}
\label{sec:introduction}

Let $\left(\Omega,\mA,\mathrm{P}\right)$ be a complete  probability space
and let
$X_i=\big(X_{1,i},\ldots X_{d,i}\big),$ $ i\geq 1,$ be  the sequence of $\bR^{d}$-valued {\it i.i.d.}  random variables  defined on $\left(\Omega,\mA, \mathrm{P}\right)$ and  having the density $f$ with respect to lebesgue measure. Furthermore,
$\bP^{(n)}_f$ denotes the probability law   of  $X^{(n)}=\big(X_1,\ldots,X_n\big),\;n\in\bN^*$ and $\bE^{(n)}_f$ is the mathematical expectation with respect to $\bP_f$.

The objective   is to estimate the density $f$ and the quality of any estimation  procedure, i.e. $X^{(n)}$-measurable mapping
$\widehat{f}_n:\bR^d\to\bL_1(\bR^d)$, is measured by {\it sup-norm risk} given by
$$
R_n^{(q)}\big(\widehat{f},f\big)=\left(\bE^{(n)}_f\big\|\widehat{f}_n-f\big\|^q_\infty\right)^{\frac{1}{q}},\;\; q\geq 1.
$$
It is well-known that even asymptotically ($n\to\infty$) the quality of estimation given by  $R_n^{(q)}$ heavily depends on the dimension
$d$. However, this asymptotics can be essentially improved if the underlying density possesses some special structure. Let us briefly discuss one of these possibilities which will be exploited in the sequel.

Introduce the following notations. Let $\cI_d$ be the set of all subsets of $\{1,\ldots,d\}$.
For any $\mathbf{I}\in\cI_d$  denote $x_{\mathbf{I}}=\left\{x_j\in\bR,\; j\in\mathbf{I}\right\}$, $\;\bar{\mathbf{I}}=\{1,\ldots,d\}\setminus\mathbf{I}$ and let $|\mathbf{I}|=\text{card}(\mathbf{I})$. Moreover for any function
$g:\bR^{|\mathbf{I}|}\to\bR$ we denote $\|g\|_{\mfI,\infty}=\sup_{x_{\mathbf{I}}\in\bR^{|\mathbf{I}|}}|g(x_{\mathbf{I}})|$.
Define also
$$
f_{\mathbf{I}}\big(x_{\mathbf{I}}\big)=\int_{\bR^{|\bar{\mathbf{I}}|}}f(x)\rd x_{\bar{\mfI}},\;\;x_{\mathbf{I}}\in\bR^{|\mathbf{I}|}.
$$
In accordance with this definition we put $f_{\mathbf{I}}\equiv 1,\;\; \mathbf{I}=\emptyset$. As  we see $f_{\mathbf{I}}$ is the marginal density of $X_{\mathbf{I},1}:=\left\{X_{j,1},\;j\in\mathbf{I}\right\}$.
Denote by $\mP$ the set of all partitions of $\{1,\ldots,d\}$  completed by  empty set $\emptyset$ and we will use $\bar{\emptyset}$
for $\{1,\ldots,d\}$.
For any density  $f$ let
$$
\mP(f)=\bigg\{\cP\in\mP:\;\; f(x)=\prod_{\mfI\in\cP}f_{\mfI}(x_\mfI),\;\;\forall x\in\bR^d\bigg\}.
$$
First we note that  $f\equiv f_{\bar{\emptyset}}$ and, therefore $\mP(f)\neq\emptyset$  since $\bar{\emptyset}\in\mP(f)$ for any $f$.
Next, if $\cP\in\mP(f)$ then $\{X_{\mfI,1},\mfI\in\cP\}$ are independent random vectors. At last, if $X_{1,1},\ldots X_{d,1},$ are independent
random variables then obviously $\mP(f)=\mP$.

Suppose now that there exists $\blP\neq\bar{\emptyset}$ such that $\blP\in\mP(f)$. If this partition is known we can proceed as follows.
For any $\blI\in\blP$ basing on observation $X^{(n)}_{\blI}$ we estimate first the marginal density$f_{\mfI}$ by $\widehat{f}_{\mfI,n}$ and then construct the estimator for
joint density $f$ as
$$
\widehat{f}_n(x)=\prod_{\blI\in\blP}\widehat{f}_{\blI,n}\left(x_{\blI}\right).
$$
One can expect (and we will see that our conjecture is true) that quality of estimation provided by this estimator will correspond not to the dimension
$d$ but to so-called effective dimension, which in our case is defined as $d(\blP)=\sup_{\blI\in\blP}|\blI|$. The main difficulty we meet trying to realize the latter construction is that the knowledge of $\blP$ is not available. Moreover, our structural hypothesis cannot be true in general, that is expressed formally by $\mP(f)=\big\{\bar{\emptyset}\big\}$. So, one of the problem we address in the present paper consists in adaptation to unknown configuration $\blP\in\mP(f)$.

We note however that even if $\blP$ is known, for instance, $\blP=\bar{\emptyset}$ the quality of an estimation procedure depends often on approximation properties of $f$ or $\{\widehat{f}_{\blI,n},\;\blI\in\blP\}$. So, our second goal is to construct an estimator which would mimic
an estimator corresponding to the minimal, and therefore unknown, approximation error. Using modern statistical language our goal here is to mimic an oracle. It is important to emphasize that we would like to solve both aforementioned problem simultaneously. Let us  now proceed with detailed consideration.

\paragraph{Collection of estimators}

Let $\mathbf{K}:\bR\to\bR$ be a given function satisfying the following assumption.
\begin{assumption}
\label{ass:on-kernel}
 $\int \mathbf{K}=1,\;$ $\|\mathbf{K}\|_\infty<\infty\;$, $\text{supp}(\mathbf{K})\subseteq [-1/2,1/2]\;$, $\mathbf{K}$ is symmetric,
 and
$$
\exists L>0:\quad  \left|\mathbf{K}(t)-\mathbf{K}(s)\right|\leq L|t-s|,\quad\forall t,s\in\bR.
$$
\end{assumption}

Put  for $\mfI\in\cI_d$
$$
K_{h_{\mfI}}(u)=V^{-1}_{h_{\mfI}}\prod_{j\in\mfI}\mathbf{K}\big(u_j/h_j\big),\qquad V_{h_{\mfI}}=\prod_{j\in\mathbf{I}}h_j.
$$
For two vectors $u,v$ here and later  $u/v$
 denotes coordinate-vise division. We will  use the notation $V_h=\prod_{j=1}^dh_j$ instead of $V_{h_{\mfI}}$ then $\mfI=\{1,\ldots,d\}$.
 Denote also $\mathrm{k}_m=\|\mathbf{K}\|_m, \;m=\{1,\infty\}$.

\smallskip

For any $p\geq 1$ let  $\gamma_p:\bN^*\times\bR_+\to\bR_+$  be the function whose explicit expression is
given in Section \ref{sec:quantity-gamma} (it has quite cumbersome expression and  it is not  convenient for us to present it right now).

Introduce the  notations (remind that $q$ is the quantity involved in the definition of the risk)
$$
\cH_n=\big\{h\in(0,1]^d:\; nV_h \geq (\ma^*)^{-1}\ln(n)\big\},\qquad \ma^*=\inf_{\mfI\in\cI_d}\left[2\gamma_{2q}\big(|\mfI|,\mathrm{k}_\infty\big)\right]^{-2}
$$
and for any $\mathbf{I}\in\cI_d$ and $h\in\cH_n$ consider kernel estimator
$$
\widetilde{f}_{h_\mathbf{I}}\big(x_{\mathbf{I}}\big)=n^{-1}\sum_{i=1}^nK_{h_{\mathbf{I}}}\left(X_{\mathbf{I},i}-x_{\mathbf{I}}\right).
$$
Introduce the family of estimators
$$
\mathfrak{F}(\mathfrak{P})=
\bigg\{\widehat{f}_{h,\cP}(x)=\prod_{\mathbf{I}\in\cP}
\widetilde{f}_{h_\mathbf{I}}\big(x_{\mathbf{I}}\big),\;\;x\in\bR^d,\;\;\cP\in\mP,\;h\in\cH_n\bigg\}.
$$
In particular,
 $
 \widehat{f}_{h,\bar{\emptyset}}(x)=n^{-1}\sum_{i=1}^nK_{h}\left(X_{i}-x\right),\;x\in\bR^d,
 $
 is the Parzen-Rosenblatt estimator (\cite{rosenblatt,parzen})   with kernel $K$ and multi-bandwidth $h$.
Our goal is to propose a data-driven selection from the family $\mF(\mP)$.

\smallskip

The estimation of a probability density is the subject of the vast literature. We do not pretend here to provide with
complete overview and only present  the results relevant in  context of the considered problems.
Minimax and minimax adaptive  density estimation
with $\bL_s$--risks   was considered in
\cite{bretagnolle}, \cite{Ibr-Has1},
\cite{devroye-gyorfi}, \cite{Efr,Efr2}, \cite{Has-Ibr},
\cite{Donoho}, \cite{Gol}, \cite{kerk},
\cite{Juditsky}, \cite{rigollet}, \cite{mason}, \cite{RRT} and \cite{akakpo}, where further references
can be found.
Oracle
inequalities for $\bL_s$--risks for $s=1$ and $s=2$ were
established  in
\cite{dev-lug96},
\cite{massart}[Chapter~7],
\cite{samarov},
\cite{rigollet-tsybakov} and \cite{birge}. The last cited paper
contains a detailed discussion of recent developments in this area. Bandwidth selection problem in the density estimation on $\bR^d$ with $\bL_s$--risks for any $1\leq s<\infty$ was studied in \cite{GL3}. The oracle inequalities obtained there were used for deriving adaptive minimax results over the collection of anisotropic Nikolskii classes.

The adaptive estimation under sup-norm loss was initiated in \cite{Lep1991,Lep1992} and continued in \cite{Tsyb} in the framework of gaussian white noise model.
Then, it was developed for anisotropic functional classes in \cite{bertin2}. The adaptive estimation of a probability density on $\bR$ in sup-norm  was the subject of recent papers \cite{GN-1,GN-2}.

\paragraph{Organization of the paper} In Section \ref{sec:oracle-supnorm} we present data-driven selection  procedure from $\mF(\mP)$ and establish for it
sup-norm oracle inequality. Section \ref{sec:supnorm-adaptation} is devoted to the adaptive estimation over the collection of anisotropic Nikolskii classes of functions. The proof of main results are given in Section \ref{sec:supnorm-main-results} and technical lemmas are proven in
Appendix.

\section{Oracle inequality}
\label{sec:oracle-supnorm}

Let $\cP\in\mP$ be fixed and define for any  $h,\eta\in\cH_n$ and any $\mfI\in\cP$
$$
\widetilde{f}_{h_\mathbf{I},\eta_{\mfI}}\big(x_{\mathbf{I}}\big)=n^{-1}\sum_{i=1}^n\big[K_{h_{\mathbf{I}}}\star K_{\eta_{\mathbf{I}}}\big] \left(X_{\mathbf{I},i}-x_{\mathbf{I}}\right),
$$
where $\big[K_{h_{\mathbf{I}}}\star K_{\eta_{\mathbf{I}}}\big]=\prod_{j\in\mfI}\big[\mathbf{K}_{h_j}\ast\mathbf{K}_{\eta_j}\big]$ and  
$
\big[\mathbf{K}_{h_j}\ast\mathbf{K}_{\eta_j}\big](z)=\int_{\bR}\mathbf{K}_{h_j}(u-z)\mathbf{K}_{\eta_j}(u)\rd u, \;\; z\in\bR.
$

\noindent As we see $"\star"$ is the convolution operator on $\bR^{|\mfI|}$.
Define
\begin{eqnarray*}
&&\mathbf{f}_n=\sup_{h\in\cH_n}\sup_{\mfI\in\cI_d}
\Big\|n^{-1}\sum_{i=1}^n\big|K_{h_{\mathbf{I}}}\left(X_{\mathbf{I},i}-\cdot\right)\big|\Big\|_{\mfI,\infty},\quad \bar{\mathbf{f}}_n=1\vee 2\mathbf{f}_n
\\*[1mm]
&&\widehat{A}_{n}(h,\cP)=\sqrt{\frac{\bar{\mathbf{f}}_n\ln(n)}{n V(h,\cP)}},\quad V(h,\cP)=\inf_{\mfI\in\cP}V_{h_{\mfI}}.
\end{eqnarray*}
Let us endow the set $\mP$ with the operation "$\diamond$" putting for any $\cP,\cP^\prime\in\mP$ 
$$
\cP\diamond\cP^\prime =\left\{\mfI\cap\mfI^\prime\neq\emptyset,\;\;\mfI\in\cP,\;\mfI^\prime\in\cP^\prime\right\}\in\mP.
$$
Introduce for any $h,\eta\in\cH_n$ and any $\cP,\cP^\prime$ the estimator
$$
\widehat{f}_{(h,\cP),(\eta,\cP^\prime)}(x)=\prod_{\mfI^{\diamond}\in\cP\diamond\cP^\prime}
\widetilde{f}_{h_{\mfI^\diamond},\eta_{\mfI^\diamond}}\big(x_{\mfI^\diamond}\big),\;\;x\in\bR^d.
$$
Set finally
$
\Lambda=\sup_{\cP\in\mP}\sup_{\mfI\in\cP}\gamma_{2q}\big(|\mfI|,\mathrm{k}_\infty\big)
$
and let $\lambda=\Lambda d\big(\bar{\mathbf{f}}_n\big)^{\left\lfloor d^{2}/4\right\rfloor+1}$.

\subsection{Selection procedure}

 For any $\cP\in\mP$ and $h\in\cH_n$  set
$$
\widehat{\Delta}_n(h,\cP)=\sup_{\eta\in\cH_n}\sup_{\cP^\prime\in\mP}
\bigg[\Big\|\widehat{f}_{(h,\cP),(\eta,\cP^\prime)}
-\widehat{f}_{\eta,\cP^\prime}\Big\|_{\infty}-\lambda\widehat{A}_{n}\big(\eta,\cP^{\prime}\big)
\bigg]_+,
$$
 and let $\widehat{h}$ and $\widehat{\cP}$ be defined as follows.
$$
\widehat{\Delta}_n\big(\widehat{h},\widehat{\cP}\big)+\lambda\widehat{A}_{n}\big(\widehat{h},\widehat{\cP}\big)=
\inf_{h\in\cH_n}\inf_{\cP\in\mP}\left[\widehat{\Delta}_n\big(h,\cP
\big)+\lambda\widehat{A}_{n}\big(h,\cP\big)\right].
$$
Our final estimator is $\widehat{f}_{\widehat{h},\widehat{\cP}}(x),\;\;x\in\bR^d$.

\paragraph{Existence and  measurability} Let us briefly discuss the existence of the proposed estimator as well as its the  measurability  with respect to the $\sigma$-algebra generated by
$X^{(n)}$. First, we note that all considered in the paper
random fields have continuous trajectories  on $\cH_n\times\bR^d$ in the topology generated by supremum norm. It is guaranteed by Assumption \ref{ass:on-kernel}. Since $\cH_n$ is totally bounded and $\bR^d$ can be covered by a countable collection of totally bounded sets, any supremum  over $\cH_n\times\bR^d$ of considered random fields will be $X^{(n)}$-measurable. In particular, $\bar{\mathbf{f}}_n$ and
$$
\widehat{\Delta}_n(h,\cP,\cP^\prime):=\sup_{\eta\in\cH_n}\sup_{\cP^\prime\in\mP}
\bigg[\Big\|\widehat{f}_{(h,\cP),(\eta,\cP^\prime)}
-\widehat{f}_{\eta,\cP^\prime}\Big\|_{\infty}-\lambda\widehat{A}_{n}\big(\eta,\cP^{\prime}\big)
\bigg]_+,\quad \cP,\cP^\prime\in\mP,\;h\in\cH_n.
$$
Since, $\mP$ is finite, we conclude that $\widehat{\Delta}_n(h,\cP)$ is $X^{(n)}$-measurable for any $\cP\in\mP$ and any $h\in\cH_n$. Assumption \ref{ass:on-kernel}  implies  also that $\widehat{\Delta}_n(\cdot,\cP)$ and $\widehat{A}_{n}\big(\cdot,\cP\big)$ are continuous on $\cH_n$  for any $\cP$. Since $\cH_n$ is a compact subset of
$\bR^d$ we conclude that $\widehat{h}(\cP)\in\cH_n$ and  $X^{(n)}$-measurable for any $\cP\in\mP$, \cite{jennrich}, where
$
\widehat{h}(\cP)=\inf_{h\in\cH_n}\left[\widehat{\Delta}_n\big(h,\cP
\big)+\lambda\widehat{A}_{n}\big(h,\cP\big)\right].
$
Since $\mP$ is finite we conclude that $(\widehat{h},\widehat{\cP})\in\cH_n\times\mP$  is  $X^{(n)}$-measurable.

\subsection{Main result} Let $\mathbf{f}>0$ be a given number and introduce the following   set of densities
$$
\mathbf{F}(\mathbf{f})=\left\{f:\;\; \sup_{\mfI\in\cI_d}\|f_{\mfI}\|_\infty\leq \mathbf{f}\right\}.
$$
With any density $f\in\mathbf{F}(\mathbf{f})$, any $h\in (0,1]^{d}$ and $\mfI\in\cI_d$ associate the quantity
$$
b_{h_{\mfI}}:=\bigg\|\int_{\bR^{|\mfI|}}K_{h_{\mfI}}\big(t_\mfI-\cdot\big)\big[f_{\mathbf{I}}
\big(t_{\mathbf{I}}\big)-f_{\mathbf{I}}(\cdot)\big]\rd t_{\mfI}\bigg\|_{\mfI,\infty},
$$
which can be view as the approximation error of  $f_{\mathbf{I}}$.

For any  $h\in\cH_n$ and $\cP\in\mP$ set
$
B\big(h,\cP\big)=\displaystyle{\sup_{\cP^\prime}\sup_{\mfI\in\cP\diamond\cP^\prime}}\left\|b_{h_{\mfI}}\right\|_{\mfI,\infty}
$
and introduce the quantity
$$
\mathfrak{R}_n(f)=\inf_{h\in\cH_n}\inf_{\cP\in\mP(f)}\left(B\big(h,\cP\big)+\sqrt{\frac{\ln(n)}{n V(h,\cP)}}\;\right).
$$

\begin{theorem}
\label{th:sup-norm-oracle}
Let Assumption \ref{ass:on-kernel} be  fulfilled. Then for any $q\geq 1$ and any $0<\mathbf{f}<\infty$ there exist
$\mathbf{C_1}\big(q,d,\mathbf{K},\mathbf{f}\big)$ and $\mathbf{C_2}\big(q,d,\mathbf{K},\mathbf{f}\big)$ such that for any $f\in\mathbf{F}(\mathbf{f})$ and any $n\geq 3$
$$
\left(\bE_f\big\|\widehat{f}_{\widehat{h},\widehat{\cP}}-f\big\|^q_\infty\right)^{\frac{1}{q}}\leq
\mathbf{C_1}\big(q,d,\mathbf{K},\mathbf{f}\big)\mathfrak{R}_n(f)+\mathbf{C_2}\big(q,d,\mathbf{K},\mathbf{f}\big)n^{-1/2}.
$$

\end{theorem}

The explicit expression of $\mathbf{C_1}\big(q,d,\mathbf{K},\mathbf{f}\big)$ and $\mathbf{C_2}\big(q,d,\mathbf{K},\mathbf{f}\big)$ can be found in the proof of the theorem.

\paragraph{Discussion}
Let us briefly discuss the assertion of Theorem \ref{th:sup-norm-oracle}.
We start with the following simple observation. Let $\bar{\mP}$ be an arbitrary subset of $\mP$ containing $\bar{\emptyset}$. If our selection rule run $\bar{\mP}$ instead of $\mP$
then the result of the  theorem remains valid if one replaces the quantity
$\mathfrak{R}_n(f)$ by
$$
\bar{\mathfrak{R}}(f)=\inf_{h\in\cH_n}\inf_{\cP\in\bar{\mP}(f)}\left(B\big(h,\cP\big)+\sqrt{\frac{\ln(n)}{n V(h,\cP)}}\;\right),
$$
where $\bar{\mP}(f)=\mP(f)\cap\bar{\mP}$. The reason of considering $\bar{\mP}$ instead of $\mP$ is explained by the fact that the cardinality of $\mP$
 (Bell number) grows as $(d/\ln(d))^d$. Therefore, for large dimension our procedure is not practically feasible in view of huge amount  of comparisons  to be done. On the other hand if $d$ is large the consideration of all partitions is not  reasonable. Indeed, even theoretically the best attainable  trade-off between approximation and stochastic errors  corresponds to the effective dimension defined as
 $
 d^{*}(f)=\inf_{\cP\in\mP(f)}\sup_{\mfI\in\cP}|\mfI|.
 $
Of course $d^{*}(f)\leq d$ but if it is proportional for example to $d$ then we will not win much for reasonable sample size. The suitable strategy in the case of large dimension consists in considering only partitions satisfying
 $
 \sup_{\mfI\in\cP}|\mfI|\leq d_0,
 $
where $d_0$ is chosen in accordance with $d$ and the number of observation. In particular one can consider $\bar{\mP}$ containing only 2 elements namely $\bar{\emptyset}$ and $\big(\{1\},\{2\},\ldots \{d\}\big)$. It corresponds to the hypotheses that we observe vectors with independent components.

Of course the consideration of $\bar{\mP}$ instead of $\mP$ has a price to pay. It is possible that $\mP(f)\cap\bar{\mP}=\bar{\emptyset}$ although $\mP(f)$ contains the elements besides $\bar{\emptyset}$. However even in this case, where structural hypothesis fails or is not taken into account ($\bar{\mP}=\{\bar{\emptyset}\}$), our estimator solves completely the bandwidths selection problem in multivariate density model under sup-norm loss.

We finish this discussion with the following remark concerning the proof of Theorem \ref{th:sup-norm-oracle}.

\begin{remark}
Our selection rule is based on computation of upper functions for some special type of random processes and   the main ingredient of the proof of Theorem \ref{th:sup-norm-oracle} is exponential inequality related  to them. Corresponding results may have an independent interest and Section \ref{sec:Upper functions for kernel estimation process}  is devoted to this topic. In particular the function $\gamma_p$  involved in the construction of our selection rule
and which we present below
comes from this consideration.
\end{remark}

\subsection{Quantity $\gamma_p$} For any $a>0$,  $p\geq 1$  and $s\in\bN^*$ introduce
\label{sec:quantity-gamma}
\begin{eqnarray*}
\gamma_p(s,a)&=&4e\sqrt{2s\tau_p(s,a)\left[a+(3L/2)(a)^{s-1}\right]}+(16e/3)\left(s\left[a+(3L/2)a^{s-1}\right]\vee 8a\right)\tau_p(s,a);
\\*[2mm]
\tau_p(s,a)&=&s\big(234s\delta_*^{-2}+6.5p+5.5\big)\ln(2)+s(2p+3)+\big[108s\delta_*^{-2}\big|\log(a)\big|+36 C_s+1\big][\ln(3)]^{-1}.
\end{eqnarray*}
Here   $\delta_{*}$ is the smallest solution of the equation
$
8\pi^2\delta\big(1+[\ln{\delta}]^{2}\big)=1
$, $C_{s}=C^{(1)}_{s}+C^{(2)}_{s}$  and
\begin{gather*}
C^{(1)}_s=s\sup_{\delta>\delta_*}\delta^{-2}\left\{\left[1+\ln{\left(\frac{9216(s+1)\delta^{2}}{[\phi(\delta)]^{2}}\right)}\right]_++1.5
\left[\log_2{\left\{\left(\frac{4608(s+1)\delta^{2}}{[\phi(\delta)]^{2}}\right)\right\}}\right]_+\right\};
\\*[2mm]
C^{(2)}_{s}=s\sup_{\delta>\delta_*}\delta^{-1}\left\{\left[1+\ln{\left(\frac{9216(s+1)\delta}{\phi(\delta)}\right)}\right]_+
+1.5\left[\log_2{\left\{\left(\frac{4608(s+1)\delta}{\phi(\delta)}\right)\right\}}\right]_+\right\},
\end{gather*}
where  $\phi(\delta)=(6/\pi^2)\big(1+[\ln{\delta}]^{2}\big)^{-1},\;\delta>0$.

\smallskip

\section{Adaptive Estimation}
\label{sec:supnorm-adaptation}

In this section we illustrate the use of the oracle inequality proved in
 Theorem \ref{th:sup-norm-oracle}  for the derivation
of adaptive rate optimal density estimators.
\par
We start with the definition of the {\em anisotropic
Nikol'skii class of functions} on $\bR^s,\;s\geq  1,$ and later on $\mathbf{e}_1,\ldots \mathbf{e}_s,$  denotes the canonical basis in $\bR^s$.
\begin{definition}
Let $r=(r_1,\ldots,r_s), r_i\in [1,\infty]$, $\alpha=(\alpha_1,\ldots, \alpha_s)$, $\alpha_i>0$, and
$Q=(Q_1,\ldots,Q_s),\;$ $Q_i>0$.
A function  $g:\bR^s\to \bR$ belongs to the
anisotropic Nikol'ski class $\bN_{r, s}(\alpha,Q)$ of functions
if
\begin{eqnarray*}
&&\|D_i^{k} g\|_{r_i} \leq  Q_i,\quad \forall k=\overline{0,\lfloor\alpha_i\rfloor},\;\;\forall i=\overline{1,s};
\\*[2mm]
&&
\left\|D_i^{\lfloor\alpha_i\rfloor} g\big(\cdot+t\mathbf{e}_i\big)- D_i^{\lfloor\alpha_i\rfloor}g\big(\cdot\big)\right\|_{r_i}\leq Q_i|t|^{\alpha_i-\lfloor\alpha_i\rfloor},\quad\forall t\in\bR,\;\;\forall i=\overline{1,s}.
\end{eqnarray*}
Here $D_i^kf$ denotes the $k$th order partial derivative of $f$ with respect to
the variable $t_i$, and $\lfloor\alpha_i\rfloor$ is the largest integer strictly
less than $\alpha_i$.
\end{definition}
\par
The functional classes $\bN_{r,s}(\alpha, Q)$ were considered in approximation
theory by
Nikol'skii; see, e.g., \cite{nikol}.
Minimax
estimation of densities from the class $\bN_{r,s}(\alpha, Q)$
was considered in
\cite{Ibr-Has2}. We refer also to
\cite{lepski-kerk,lepski-kerk2}, where the problem of adaptive estimation over a
scale of classes
$\bN_{r,s}(\alpha,Q)$ was treated for the Gaussian white noise model.

\smallskip

Our goal now is to introduce the scale of functional classes of $d$-variate probability densities taking into account the  independence structure.
It implies in particular that we will need to estimate not only the density itself but all  marginal densities as well. It  is easily seen
that if  $f\in\bN_{p, d}(\beta,\cL)$ and additionally $f$ is compactly supported  then $f_{\mfI}\in\bN_{p_\mfI,|\mfI|}\big(\beta_\mfI,\overline{\cL}_\mfI\big)$ for any $\mfI\in\cI_d$, where $\overline{\cL}=\mathrm{c}\cL$ and
$\mathrm{c}>0$ is a numerical constant. However if $\text{supp}(f)=\bR^d$ the latter assertion is not true in general. The assumption
$f\in\bN_{p, d}(\beta,\cL)$ does not even guarantee that $f_\mfI$ is bounded on $\bR^{|\mfI|}$. It explains the introduction of the following
anisotropic classes of densities.

\vskip0.1cm

Let  $p=(p_1,\ldots,p_d), p_i\in [1,\infty]$, $\beta=(\beta_1,\ldots, \beta_d)$, $\beta_i>0$,
$\cL=(\cL_1,\ldots,\cL_d),\;$ $\cL_i>0$.
\begin{definition}
A probability density $f:\bR^d\to\bR_+$ belongs to the  class $\overline{\bN}_{p,d}\big(\beta,\cL\big)$ if
$$
 f_{\mfI}\in \bN_{p_\mfI,|\mfI|}\big(\beta_\mfI,\cL_\mfI\big),\;\;
\forall\mfI\in\cI_d.
$$
\end{definition}

Introduce finally the collection of functional classes taking into account the smoothness of the underlying density and the independence structure simultaneously.

 Let $\big(\beta,p,\cP\big)\in (0,\infty)^{d}\times[1,\infty]^d\times\mP$ and $\cL\in (0,\infty)^d$ be fixed. Introduce
$$
\mathbf{N}_{p,d}\big(\beta,\cL,\cP\big)=\bigg\{f(x)\in\overline{\bN}_{p,d}\big(\beta,\cL\big):\;\; f(x)=\prod_{\mfI\in\cP}f_{\mfI}\big(x_\mfI\big),\;\;\forall x\in\bR^d\bigg\}.
$$
For any $\big(\beta,p,\cP\big)\in (0,\infty)^{d}\times[1,\infty]^d\times\mP$  define
$$
\Upsilon\big(\beta,p,\cP\big)=\inf_{\mfI\in\cP}\gamma_\mfI(\beta,p),\qquad
\gamma_{\mfI}\big(\beta,p\big)=\frac{1-\sum_{j\in\mfI}\frac{1}{\beta_jp_j}}{\sum_{j\in\mfI}\frac{1}{\beta_j}}.
$$
We will see that the quantity  $\Upsilon\big(\beta,p,\cP\big)$ can be view as "effective smoothness index" related to independence structure hypothesis and to the estimation under sup-norm loss.

\begin{theorem}
\label{th:lower-bound-in-supnorm}
For any $\big(\beta,p,\cP\big)\in (0,\infty)^{d}\times[1,\infty]^d\times\mP$ such that $\Upsilon\big(\beta,p,\cP\big)>0$ and any $\cL\in (0,\infty)^d$
$$
\liminf_{n\to\infty}\;\inf_{\widehat{f}_n}\sup_{f\in\mathbf{N}_{p,d}\big(\beta,\cL,\cP\big)}
\left(\bE^{(n)}_f\left[\varphi_n^{-1}(\beta,p,\cP\big)\big\|\widehat{f}_n-f\big\|_\infty\right]^q\right)^{\frac{1}{q}}>0,\quad\;
\varphi_n(\beta,p,\cP\big)=\left(\frac{\ln n}{n}\right)^{\frac{\Upsilon}{2\Upsilon+1}}.
$$
where $\Upsilon=\Upsilon\big(\beta,p,\cP\big)$ and infimum is taken over all possible estimators.
\end{theorem}

Our goal is to prove that the estimation quality provided by  $\widehat{f}_{\widehat{h},\widehat{\cP}}$ on
$\mathbf{N}_{p,d}\big(\beta,\cL,\cP\big)$  coincides up to numerical constant with optimal decay of minimax risk  $\varphi_n(\beta,p,\cP\big)$ whenever the value of
nuisance parameter $\big\{\beta,p,\cP, \cL\big\}$. It means that this estimator is optimally adaptive  over the scale of considered functional classes.
We would like to emphasize that not only the couple $(\beta,\cL)$ is unknown that is typical in frameworks of adaptive estimation but also the index $p$ of  norms where the smoothness is measured. At last, our estimator adapts automatically to unknown independence structure.

\begin{theorem}
\label{th:adaptation-in-supnorm}
Let $\mathbf{K}$ satisfy Assumption \ref{ass:on-kernel} and suppose additionally that for some $\mb> 2$
\begin{equation}
\label{eq:ass-th:adaptation-in-supnorm}
\int_{\bR}u^{m}\mathbf{K}(u)\rd u=0,\quad\forall m=\overline{2,\mb}.
\end{equation}

Then for any $\big(\beta,p,\cP\big)\in (0,\mb]^{d}\times[1,\infty]^d\times\mP$ such that $\Upsilon\big(\beta,p,\cP\big)>0$ and any $\cL\in (0,\infty)^d$
$$
\limsup_{n\to\infty}\sup_{f\in\mathbf{N}_{p,d}\big(\beta,\cL,\cP\big)}
\left(\bE^{(n)}_f\left[\varphi_n^{-1}(\beta,p,\cP\big)\big\|\widehat{f}_{\widehat{h},\widehat{\cP}}
-f\big\|_\infty\right]^q\right)^{\frac{1}{q}}<\infty.
$$
\end{theorem}
 We want to emphasize that the extra-parameter $\mb$ can be arbitrary but {\it a priory} chosen.
 Note  that the condition (\ref{eq:ass-th:adaptation-in-supnorm}) of the theorem is fulfilled with $m=1$ as well  since $\mathbf{K}$ is symmetric.

We remark also that  for any given $\big(\beta,p,\cP\big)\in (0,\mb]^{d}\times[1,\infty]^d\times\mP$, satisfying $\Upsilon\big(\beta,p,\cP\big)>0$,  one can find $\mathbf{f}=\mathbf{f}\big(\beta,p,\cP\big)$ such that
$f\in\mathbf{N}_{p,d}\big(\beta,\cL,\cP\big)$ implies that $f\in\mathbf{F}(\mathbf{f})$. It makes possible the application of Theorem \ref{th:sup-norm-oracle}.

\section{Proofs}
\label{sec:supnorm-main-results}

We start this section with the computation of upper functions  for kernel estimation process being one of  main tools in the proof of Theorem \ref{th:sup-norm-oracle}.

\subsection{Upper functions for kernel estimation process }
\label{sec:Upper functions for kernel estimation process}
Let $s\in\bN^*$ and let $Y_j, j\geq 1,$ be $\bR^s$-valued i.i.d. random vectors defined on a complete  probability space $\left(\Omega,\mA,\mathrm{P}\right)$
and having the density $\mathbf{g}$ with respect to the Lebesgue measure.
Later on  $\bP^{(n)}_\mathbf{g}$ denotes the law of $Y_1,\ldots, Y_n, n\in\bN^*,$ and $\bE^{(n)}_\mathbf{g}$  is mathematical expectation with respect to $\bP^{(n)}_\mathbf{g}$.

Let $\mathbf{M}:\bR\to\bR$ be a given symmetric function and for any $r\in (0,1]^s$ set as previously
$$
M_r(\cdot)=\prod_{l=1}^s r^{-1}_l\mathbf{M}(\cdot/r_l),\quad V_r=\prod_{l=1}^s r_l.
$$
Denote also $\mathrm{m}_m=\|\mathbf{M}\|_m, \;m=\{1,\infty\}$.
For any $y\in\bR^s$ consider the family of random fields
$$
\chi_r(y)=n^{-1}\sum_{j=1}^n\left\{M_{r}\left(Y_j-y\right)-\bE^{(n)}_\mathbf{g}\Big[M_r\left(Y_j-y\right)\Big]\right\},\;\; r\in\widetilde{\cR}_n(s):=\left\{r\in (0,1]^s:\;\; nV_r\geq \ln(n)\right\}.
$$
For  any $r\in (0,1]^s$ set
$
G(r)=\displaystyle{\sup_{y\in\bR^s}\int_{\bR^s}}|M_r(x-y)|\mathbf{g}(x)\rd x
$
and let $\bar{G}(r)=1\vee G(r)$.

\begin{proposition}
\label{prop:exp-ineq-for-kernel-process}
Let $\mathbf{M}$ satisfy Assumption \ref{ass:on-kernel}. Then for any $n\geq 3$ and any $p\geq 1$
$$
\bE^{(n)}_\mathbf{g}\Bigg\{\sup_{r\in\widetilde{\cR}_n(s)}\bigg[\big\|\chi_r\big\|_\infty-\gamma_p\big(s,\mathrm{m}_\infty\big)
\sqrt{\frac{\bar{G}(r)\ln(n)}{nV_r}}\bigg]\Bigg\}^p_+\leq c_1(p,s)\big[1\vee\mathrm{m}^{s}_1\|\mathbf{g}\|_\infty\big]^{\frac{p}{2}}n^{-\frac{p}{2}}+c_2(p,s)n^{-p},
$$
where  $c_1(p,s)=2^{7p/2+5}3^{p+5s+4}\Gamma(p+1)\pi^{p}\big(s,\mathrm{m}_\infty\big)$ and $c_2(p,s)=2^{p+1}3^{5s}$.
\end{proposition}
The function $\pi:\bN^*\times\bR_+:\to\bR$ is given by
$$
\pi(s,a)=\big(\sqrt{a}\vee a\big)\left(\sqrt{2es\left[1+(3L/2)a^{s-2}\right]}\vee \bigg[(2e/3)\bigg(s\left[1+(3L/2)a^{s-2}\right]\vee 8\bigg)\bigg]\right).
\\
$$

In view of trivial inequality
$$
\big\|\chi_r\big\|_\infty\leq \gamma_p\big(s,\mathrm{m}_\infty\big)
\sqrt{\frac{\bar{G}(r)\ln(n)}{nV_r}}+\left(\big\|\chi_r\big\|_\infty-\gamma_p\big(s,\mathrm{m}_\infty\big)
\sqrt{\frac{\bar{G}(r)\ln(n)}{nV_r}}\right)_+
$$
we come to the following corollary of Proposition \ref{prop:exp-ineq-for-kernel-process}.
\begin{corollary}
\label{cor:prop:exp-ineq-for-kernel-process}
Let $\mathbf{M}$ satisfy Assumption \ref{ass:on-kernel}. Then for any $n\geq 3$ and any $p\geq 1$
$$
\Bigg(\bE^{(n)}_\mathbf{g}\bigg\{\sup_{r\in\widetilde{\cR}_n(s)}\big\|\chi_r\big\|_\infty\bigg\}^p\Bigg)^{\frac{1}{p}}\leq
\big[1\vee\mathrm{m}^{s}_1\|\mathbf{g}\|_\infty\big]^{\frac{1}{2}}
\left[\gamma_p\big(s,\mathrm{m}_\infty\big)+\big\{c_1(p,s)+c_2(p,s)\big\}^{\frac{1}{p}}n^{-1/2}\right].
$$

\end{corollary}

Consider now the following family of random processes: for any $y\in\bR^s$
$$
\Upsilon_r(y)=n^{-1}\sum_{j=1}^n\big|M_{r}\left(Y_j-y\right)\big|,\;\; r\in\widetilde{\cR}^{(\ma)}_n(s):=\left\{r\in (0,1]^s:\;\; nV_r\geq \ma^{-1}\ln(n)\right\},
$$
where we have put $\ma=\left[2\gamma_p\big(s,\mathrm{m}_\infty\big)\right]^{-2}$.

\begin{proposition}
\label{prop2:exp-ineq-for-kernel-process}
Let $\mathbf{M}$ satisfy Assumption \ref{ass:on-kernel}. Then for any $n\geq 3$ and any $p\geq 1$
\begin{eqnarray*}
&&\bE^{(n)}_\mathbf{g}\bigg\{\sup_{r\in\cR_n^{(\ma)}(s)}\left[1\vee\left\|\Upsilon_r\right\|_\infty-(3/2)\bar{G}(r)\right]\bigg\}^p_+\leq
c_1(p,s)\big[1\vee\mathrm{m}^{s}_1\|\mathbf{g}\|_\infty\big]^{\frac{p}{2}}n^{-\frac{p}{2}}+c_2(p,s)n^{-p};
\\
&&\bE^{(n)}_\mathbf{g}\bigg\{\sup_{r\in\cR_n^{(\ma)}(s)}\left[\bar{G}(r)-
2\left(1\vee\left\|\Upsilon_r(\cdot)\right\|_\infty\right)\right]\bigg\}^p_+\leq c^{\prime}_1(p,s)\big[1\vee\mathrm{m}^{s}_1\|\mathbf{g}\|_\infty\big]^{\frac{p}{2}}n^{-\frac{p}{2}}+c^{\prime}_2(p,s)n^{-p},
\end{eqnarray*}
where  $c^{\prime}_1(p,s)=2^pc_1(p,s)$ and $c^{\prime}_2(p,s)=2^{2p+1}3^{5s}$.
\end{proposition}

\subsection{Proof of Theorem \ref{th:sup-norm-oracle}} We start the proof of the theorem
with  auxiliary results used in the sequel. Whose proofs are given in Appendix.

\subsubsection{Auxiliary results}

Introduce the following notations. For any $\mfI\in\cI_d$ set
$$
s_{h_\mfI}(\cdot)=\int_{\bR^{|\mfI|}}K_{h_{\mfI}}\big(t_\mfI-\cdot\big)f_{\mathbf{I}}
\big(t_{\mathbf{I}}\big)\rd t_\mfI, \qquad s^*_{h_\mfI,\eta_\mfI}(\cdot)=\int_{\bR^{|\mfI|}}\left[K_{h_{\mfI}}\star K_{\eta_{\mfI}}\right]\big(t_\mfI-\cdot\big)f_{\mathbf{I}}
\big(t_{\mathbf{I}}\big)\rd t_\mfI;
$$

\begin{lemma}
\label{lem:bais-supnorm-result} For any  $\mfI\in\cI_d$ and any $h,\eta\in (0,1]^{|\mfI|}$ one has
$$
\left\|s^*_{h_\mfI,\eta_\mfI}-s_{\eta_\mfI}\right\|_{\mfI,\infty}\leq \mathrm{k}^{d}_1b_{h_{\mfI}}.
$$

\end{lemma}


For any $h\in (0,1]^{d}$ and any $\cP\in\mP$ set
$$
A_{n}(h,\cP)=\sqrt{\frac{\bar{s}_n\ln(n)}{n V(h,\cP)}},\quad\bar{s}_n=1\vee\sup_{h\in\cH_n }\sup_{\mfI\in\cI_d}\bigg\|\int_{\bR^{|\mfI|}}\big|K_{h_{\mfI}}\big(t_\mfI-\cdot\big)\big|f_{\mathbf{I}}
\big(t_{\mathbf{I}}\big)\rd t_\mfI\bigg\|_{\mfI,\infty}
$$
Put also   $\xi_{h_\mfI}(\cdot)=\widetilde{f}_{h_\mathbf{I}}(\cdot)-s_{h_\mfI}(\cdot)$ and let
$$
\zeta(h,\cP)
=\sup_{\mfI\in\cP}\left\|\xi_{h_\mfI}\right\|_{\mfI,\infty},
\quad\zeta_n=\sup_{\eta\in\cH_n}\sup_{\cP\in\mP}\bigg[\zeta\big(\eta,\cP\big)-\Lambda A_n(\eta,\cP)\bigg]_+,
$$
\begin{lemma}
\label{lem:stoch-part}

For any $p\geq 1$ there exist $\mathbf{c_i}\big(p,d,\mathbf{K},\mathbf{f}\big),\;\mfi=1,2,3,4,$
such that for any $n\geq 3$
\begin{eqnarray*}
&&(\mfi)\qquad\sup_{f\in\mathbf{F}(\mathbf{f})}\left[\bE^{(n)}_f\big(\zeta_n\big)^{2q}\right]^{\frac{1}{2q}}\leq \mathbf{c_1}\big(2q,d,\mathbf{K},\mathbf{f}\big)n^{-1/2};
\\
&&(\mfi\mfi)\qquad\sup_{f\in\mathbf{F}(\mathbf{f})}\left[\bE^{(n)}_f\left[\bar{s}_n-\bar{\mathbf{f}}_n\right]_+^{2q}\right]^{\frac{1}{2q}}\leq \mathbf{c_2}\big(2q,d,\mathbf{K},\mathbf{f}\big)n^{-1/2};
\\
&&(\mfi\mfi\mfi)\qquad\sup_{f\in\mathbf{F}(\mathbf{f})}\left[\bE^{(n)}_f\left[\bar{\mathbf{f}}_n-3\bar{s}_n\right]_+^{2q}\right]^{\frac{1}{2q}}\leq \mathbf{c_3}\big(2q,d,\mathbf{K},\mathbf{f}\big)n^{-1/2};
\\
&&(\mfi\mathbf{v})\qquad\sup_{f\in\mathbf{F}(\mathbf{f})}\left[\bE^{(n)}_f\big(\bar{\mathfrak{f}}_n\big)^p\right]^{\frac{1}{p}}\leq \mathbf{c_4}\big(p,d,\mathbf{K},\mathbf{f}\big).
\end{eqnarray*}

\end{lemma}

The explicit expression of $\mathbf{c_i}\big(p,d,\mathbf{K},\mathbf{f}\big),\;\mfi=1,2,3,4$ can be found in the proof of the lemma.

\subsubsection{Proof of Theorem \ref{th:sup-norm-oracle}}
We brake the proof on several steps.

\smallskip

$\mathbf{1^0}.$ Let  $\blh\in\cH_n$ and $\blP\in\mP$ be fixed.
We have in view of triangle inequality
\begin{equation}
\label{eq02:proof-th:sup-norm-oracle}
\Big\|\widehat{f}_{\widehat{h},\widehat{\cP}}-f\Big\|_\infty\leq
\Big\|\widehat{f}_{\widehat{h},\widehat{\cP}}-\widehat{f}_{(\blh,\blP)(\widehat{h},\widehat{\cP})}\Big\|_\infty+
\Big\|\widehat{f}_{(\blh,\blP)(\widehat{h},\widehat{\cP})}-\widehat{f}_{\blh,\blP}\Big\|_\infty+
\Big\|\widehat{f}_{\blh,\blP}-f\Big\|_\infty.
\end{equation}
We have
\begin{equation}
\label{eq03:proof-th:sup-norm-oracle}
\Big\|\widehat{f}_{\widehat{h},\widehat{\cP}}-\widehat{f}_{(\blh,\blP)(\widehat{h},\widehat{\cP})}\Big\|_\infty
\leq \widehat{\Delta}_n(\blh,\blP)+\lambda\widehat{A}_{n}\big(\widehat{h},\widehat{\cP}\big).
\end{equation}
Noting that $\widehat{f}_{(\blh,\blP)(\widehat{h},\widehat{\cP})}\equiv \widehat{f}_{(\widehat{h},\widehat{\cP})(\blh,\blP)}$ we get
\begin{equation}
\label{eq05:proof-th:sup-norm-oracle}
\Big\|\widehat{f}_{(\blh,\blP)(\widehat{h},\widehat{\cP})}-\widehat{f}_{\blh,\blP}\Big\|_\infty\leq \widehat{\Delta}_n\big(\widehat{h},\widehat{\cP}\big)+\lambda\widehat{A}_{n}\big(\blh,\blP\big).
\end{equation}
We obtain from  (\ref{eq03:proof-th:sup-norm-oracle}) and (\ref{eq05:proof-th:sup-norm-oracle})
\begin{eqnarray*}
\label{eq06:proof-th:sup-norm-oracle}
&&\Big\|\widehat{f}_{\widehat{h},\widehat{\cP}}-\widehat{f}_{(\blh,\blP)(\widehat{h},\widehat{\cP})}\Big\|_\infty+
\Big\|\widehat{f}_{(\blh,\blP)(\widehat{h},\widehat{\cP})}-\widehat{f}_{\blh,\blP}\Big\|_\infty
\\*[2mm]
&&\leq
\left[\widehat{\Delta}_n\big(\widehat{h},\widehat{\cP}\big)+\lambda\widehat{A}_{n}\big(\widehat{h},\widehat{\cP}\big)\right]+
\left[\widehat{\Delta}_n(\blh,\blP)+\lambda\widehat{A}_{n}(\blh,\blP)\right]\leq 2\left[\widehat{\Delta}_n(\blh,\blP)+\lambda\widehat{A}_{n}(\blh,\blP)\right].
\end{eqnarray*}
To get the last inequality we have used the definition of $(\widehat{h},\widehat{\cP})$.
Thus, we obtain from (\ref{eq02:proof-th:sup-norm-oracle}) that
\begin{equation}
\label{eq1:proof-th:sup-norm-oracle}
\Big\|\widehat{f}_{\widehat{h},\widehat{\cP}}-f\Big\|_\infty\leq 2\left[\widehat{\Delta}_n(\blh,\blP)+\lambda\widehat{A}_{n}(\blh,\blP)\right]
+\Big\|\widehat{f}_{\blh,\blP}-f\Big\|_\infty.
\end{equation}

\smallskip

$\mathbf{2^0}.$ 
Note that for any $h,\eta\in\cH_n$ and any $\cP^\prime\in\mP$
\begin{equation}
\label{eq2:proof-th:sup-norm-oracle}
\Big\|\widehat{f}_{(h,\blP),(\eta,\cP^\prime)}
-\widehat{f}_{\eta,\cP^\prime}\Big\|_{\infty}\leq
d\big(\bar{\mathbf{f}}_n\big)^{\left\lfloor d^{2}/4\right\rfloor+1}
\sup_{\mfI^{\prime}\in\cP^\prime}\bigg\|\prod_{\blI\in\blP:\;\blI\cap\mfI^\prime\neq\emptyset}
\widetilde{f}_{h_{\blI\cap\mfI^\prime},\eta_{\blI\cap\mfI^\prime}}
-\widetilde{f}_{\eta_{\mfI^\prime}}\bigg\|_{\mfI^\prime,\infty}.
\end{equation}
Here we have used the trivial inequality: for any $m\in\bN^*$ and any $a_j,b_j:\cX_j\to\bR,\;j=\overline{1,m},$
\begin{equation}
\label{eq3:proof-th:sup-norm-oracle}
\bigg\|\prod_{j=1}^m a_j-\prod_{j=1}^m b_j\bigg\|_\infty\leq m\bigg(\sup_{j=\overline{1,m}}\|a_j-b_j\|_{\cX_j,\infty}\bigg)
\bigg[\sup_{j=\overline{1,m}}\max\big(\|a_j\|_{\cX_j,\infty},\|b_j\|_{\cX_j,\infty}\big)\bigg]^{m-1},
\end{equation}
where $\|\cdot\|_{\cX_j,\infty}$ and $\|\cdot\|_{\infty}$ denote the supremum norms on $\cX_j$ and $\cX_1\times\cdots\times\cX_m$ respectively.

Introduce the following notation: for any $h,\eta\in\cH_n$ and any $\cP\in\mP$ we set
\begin{gather*}
\xi_{h_\mfI,\eta_\mfI}^*(\cdot)=\widetilde{f}_{h_\mathbf{I},\eta_{\mfI}}(\cdot)-s^*_{h_\mfI,\eta_\mfI}(\cdot)
\end{gather*}
 We have in view of (\ref{eq3:proof-th:sup-norm-oracle}) (here and later  the product and the supremum over empty set are assumed equal to one and to zero respectively)
\begin{equation}
\label{eq4:proof-th:sup-norm-oracle}
\bigg\|\prod_{\blI\in\blP}
\widetilde{f}_{h_{\blI\cap\mfI^\prime},\eta_{\blI\cap\mfI^\prime}}
-\prod_{\blI\in\blP}s^*_{h_{\blI\cap\mfI^\prime},\eta_{\blI\cap\mfI^\prime}}\bigg\|_{\mfI^\prime,\infty}\leq
d\left[\max\big\{\bar{\mathbf{f}}_n,\rk^2_1\mathbf{f}\big\}\right]^{d-1}
\sup_{\blI\in\blP}\left\|\xi_{h_{\blI\cap\mfI^\prime},\eta_{\blI\cap\mfI^\prime}}^*\right\|_{\blI\cap\mfI^\prime,\infty}.
\end{equation}
We remark that for any $\mfI\in\cI_d$, any $h,\eta\in (0,1]^d$ and any $z_{\mfI}\in\bR^{|\mfI|}$
$$
\xi_{h_{\mfI},\eta_{\mfI}}^*\big(z_{\mfI}\big)=\int_{\bR^{|\mfI|}}K_{\eta_i}\big(z_{\mfI}-u_{\mfI}\big)\xi_{h_{\mfI}}\big(u_{\mfI}\big)\rd u_{\mfI}
$$
and, therefore,
$$
\left\|\xi_{h_{\mfI},\eta_{\mfI}}^*\right\|_{\mfI,\infty}\leq \mathrm{k_1}^{|\mfI|}\left\|\xi_{h_{\mfI}}\right\|_{\mfI,\infty}\leq
 \mathrm{k_1}^{d}\left\|\xi_{h_{\mfI}}\right\|_{\mfI,\infty},
$$
since $\rk_1\geq 1$ in view of Assumption \ref{ass:on-kernel}. It yields together with (\ref{eq4:proof-th:sup-norm-oracle})
\begin{equation}
\label{eq40:proof-th:sup-norm-oracle}
\bigg\|\prod_{\blI\in\blP}
\widetilde{f}_{h_{\blI\cap\mfI^\prime},\eta_{\blI\cap\mfI^\prime}}
-\prod_{\blI\in\blP}s^*_{h_{\blI\cap\mfI^\prime},\eta_{\blI\cap\mfI^\prime}}\bigg\|_{\mfI^\prime,\infty}\leq
d \mathrm{k_1}^{d}\left[\max\big\{\bar{\mathbf{f}}_n,\rk^2_1\mathbf{f}\big\}\right]^{d-1}
\sup_{\blI\in\blP}\left\|\xi_{h_{\mfI\cap\mfI^\prime}}\right\|_{\blI\cap\mfI^\prime,\infty}.
\end{equation}

Note also that for any $\eta\in\cH_n$ and $\mfI^\prime\in\cI_d$
\begin{eqnarray*}
s_{\eta_{\mfI^\prime}}(\cdot)=\int_{\bR^{\mfI^\prime}}K_{\eta_{\mfI^\prime}}\big(t_{\mfI^\prime}-\cdot\big)f_{\mfI^\prime}
\big(t_{\mfI^\prime}\big)\rd t_{\mfI^\prime}=\int_{\bR^{\mfI^\prime}}K_{\eta_{\mfI^\prime}}\big(t_{\mfI^\prime}-\cdot\big)\bigg[\prod_{\blI\in\blP}f_{\blI\cap\mfI^\prime}
\big(t_{\blI\cap\mfI^\prime}\big)\bigg]\rd t_{\mfI^\prime}=\prod_{\blI\in\blP}s_{\eta_{\blI\cap\mfI^\prime}}(\cdot).
\end{eqnarray*}
Here we have used that $\blP\in\mP(f)$. Using once again (\ref{eq3:proof-th:sup-norm-oracle}) we obtain
\begin{equation*}
\label{eq5:proof-th:sup-norm-oracle}
\bigg\|
\prod_{\blI\in\blP}s^*_{h_{\blI\cap\mfI^\prime},\eta_{\blI\cap\mfI^\prime}}
-\prod_{\blI\in\blP}s_{\eta_{\blI\cap\mfI^\prime}}\bigg\|_{\mfI^\prime,\infty}\leq
d\left[\rk^2_1\mathbf{f}\right]^{d-1}
\sup_{\blI\in\blP}
\left\|s^*_{h_{\blI\cap\mfI^\prime},\eta_{\blI\cap\mfI^\prime}}-s_{\eta_{\blI\cap\mfI^\prime}}\right\|_{\blI\cap\mfI^\prime,\infty}.
\end{equation*}
and, therefore, in view of Lemma \ref{lem:bais-supnorm-result}
\begin{eqnarray}
\label{eq6:proof-th:sup-norm-oracle}
\bigg\|
\prod_{\blI\in\blP}s^*_{h_{\blI\cap\mfI^\prime},\eta_{\blI\cap\mfI^\prime}}
-s_{\eta_{\mfI^\prime}}\bigg\|_{\mfI^\prime,\infty}&\leq&
d\rk^d_1\left[\rk^{2}_1\mathbf{f}\right]^{d-1}
\sup_{\blI\in\blP}
\left\|b_{h_{\blI\cap\mfI^\prime}}\right\|_{\blI\cap\mfI^\prime,\infty}.
\end{eqnarray}
Thus, we obtain from (\ref{eq40:proof-th:sup-norm-oracle}) and (\ref{eq6:proof-th:sup-norm-oracle})
\begin{equation*}
\bigg\|\prod_{\blI\in\blP:\;\blI\cap\mfI^\prime\neq\emptyset}
\widetilde{f}_{h_{\blI\cap\mfI^\prime},\eta_{\blI\cap\mfI^\prime}}
-\widetilde{f}_{\eta_{\mfI^\prime}}\bigg\|_{\mfI^\prime,\infty}\leq \mathfrak{f}_n
\bigg[\sup_{\blI\in\blP}\left\|\xi_{h_{\mfI\cap\mfI^\prime}}\right\|_{\blI\cap\mfI^\prime,\infty}
+\sup_{\blI\in\blP}
\left\|b_{h_{\blI\cap\mfI^\prime}}\right\|_{\blI\cap\mfI^\prime,\infty}
\bigg]+\left\|\xi_{\eta_{\mfI^{\prime}}}\right\|_{\mfI^{\prime},\infty},
\end{equation*}
where we have put $\mathfrak{f}_n=2d\rk_1^d\left[\max\big\{\bar{\mathbf{f}}_n,\rk^2_1\mathbf{f}\big\}\right]^{d-1}$.

Therefore, we get from (\ref{eq2:proof-th:sup-norm-oracle}) for any $h,\eta\in\cH_n$ and $\cP^\prime\in\mP$

\begin{eqnarray}
\label{eq7:proof-th:sup-norm-oracle}
&&\Big\|\widehat{f}_{(h,\blP),(\eta,\cP^\prime)}
-\widehat{f}_{\eta,\cP^\prime}\Big\|_{\infty}\leq
\bar{\mathfrak{f}}_n\bigg\{
\zeta\big(h,\blP\diamond\cP^\prime\big)
+\sup_{\mfI\in\blP\diamond\cP^\prime}
\left\|b_{h_{\mfI}}\right\|_{\mfI,\infty}\bigg\}
+\tilde{\mathfrak{f}}_n\zeta\big(\eta,\cP^\prime\big).
\end{eqnarray}
Here we have put $\tilde{\mathfrak{f}}_n=d\big(\bar{\mathbf{f}}_n\big)^{\left\lfloor d^{2}/4\right\rfloor+1}$ and $\bar{\mathfrak{f}}_n=\tilde{\mathfrak{f}}_n\mathfrak{f}_n$ 
Taking into account  that for any $h\in\cH_n$ and any $\cP,\cP^\prime\in\mP$
\begin{equation*}
\label{eq005:proof-th:sup-norm-oracle}
A_{n}\big(h,\cP\diamond\cP^\prime\big)\leq A_{n}\big(h,\cP\big)\wedge A_{n}\big(h,\cP^\prime\big),
\end{equation*}
we get from (\ref{eq7:proof-th:sup-norm-oracle})
\begin{eqnarray}
\label{eq70:proof-th:sup-norm-oracle}
&&\Big\|\widehat{f}_{(h,\blP),(\eta,\cP^\prime)}
-\widehat{f}_{\eta,\cP^\prime}\Big\|_{\infty}\leq
\bar{\mathfrak{f}}_n\bigg\{
\Lambda A_{n}\big(h,\blP\big)
+\sup_{\mfI\in\blP\diamond\cP^\prime}
\left\|b_{h_{\mfI}}\right\|_{\mfI,\infty}+\zeta_n\bigg\}
+\tilde{\mathfrak{f}}_n\zeta\big(\eta,\cP^\prime\big).
\end{eqnarray}
Remembering that $\lambda=\tilde{\mathfrak{f}}_n\Lambda$, we obtain
from (\ref{eq70:proof-th:sup-norm-oracle})
\begin{equation*}
\widehat{\Delta}_n(\blh,\blP)\leq \bar{\mathfrak{f}}_n\Big\{\Lambda A_{n}\big(\blh,\blP\big)+B\big(\blh,\blP\big)
+\zeta_n\Big\}+\tilde{\mathfrak{f}}_n\Big\{\zeta_n+\Lambda\sup_{\eta\in\cH_n}\sup_{\cP\in\cP}\left[A_n(\eta,\cP)-\widehat{A}_n(\eta,\cP)\right]_+\Big\},
\end{equation*}
where, remind $B\big(h,\cP\big)=\sup_{\cP^\prime}\sup_{\mfI\in\cP\diamond\cP^\prime}
\left\|b_{h_{\mfI}}\right\|_{\mfI,\infty}$. 

Taking into account that $\bar{\mathfrak{f}}_n\geq \tilde{\mathfrak{f}}_n$, since $\bar{\mathfrak{f}}_n\geq 1$ we finally get
\begin{eqnarray}
\label{eq8:proof-th:sup-norm-oracle}
&&\widehat{\Delta}_n(\blh,\blP)\leq \bar{\mathfrak{f}}_n\Big\{\Lambda A_{n}\big(\blh,\blP\big)+B\big(\blh,\blP\big)
+2\zeta_n+\Lambda\sup_{\eta\in\cH_n}\sup_{\cP\in\cP}\left[A_n(\eta,\cP)-\widehat{A}_n(\eta,\cP)\right]_+\Big\}.
\end{eqnarray}

Note that the definition of $\cH_n$ implies that
$$
\left[A_n(\eta,\cP)-\widehat{A}_n(\eta,\cP)\right]_+\leq \ma^*\left[\sqrt{\bar{s}_n}-\sqrt{\bar{\mathbf{f}}_n}\right]_+\leq
\ma^* \left[\bar{s}_n-\bar{\mathbf{f}}_n\right]_+,\quad \forall \eta\in\cH_n,\;\forall\cP\in\cP.
$$
To get the last inequality we have also used that by definition $\bar{\mathbf{f}}_n,\bar{s}_n\geq 1.$

Putting
$
R_n=\ma^*\Lambda \left[\bar{s}_n-\bar{\mathbf{f}}_n\right]_+
$
we obtain in view of (\ref{eq8:proof-th:sup-norm-oracle})
\begin{eqnarray}
\label{eq9:proof-th:sup-norm-oracle}
&&\widehat{\Delta}_n(\blh,\blP)\leq \bar{\mathfrak{f}}_n\Big\{\Lambda A_{n}\big(\blh,\blP\big)+B\big(\blh,\blP\big)
+2\zeta_n+R_n\Big\},
\end{eqnarray}
Note also that   the definition of $\cH_n$ implies that
\begin{eqnarray*}
&&\left[\widehat{A}_{n}(\blh,\blP)-\sqrt{3}A_{n}(\blh,\blP)\right]_+\leq \ma^*\left[\sqrt{\bar{\mathbf{f}}_n}-\sqrt{3\bar{s}_n}\right]_+\leq
\ma^* \left[\bar{\mathbf{f}}_n-3\bar{s}_n\right]_+,\quad \forall \eta\in\cH_n,\;\forall\cP\in\cP.
\end{eqnarray*}
Thus, denoting $\cR_n=\ma^*\Lambda\left[\bar{\mathbf{f}}_n-3\bar{s}_n\right]_+$ we obtain using (\ref{eq9:proof-th:sup-norm-oracle})
\begin{eqnarray}
\label{eq99:proof-th:sup-norm-oracle}
&&\widehat{\Delta}_n(\blh,\blP)+\lambda\widehat{A}_{n}\big(\blh,\blP\big)\leq \bar{\mathfrak{f}}_n\Big\{3\Lambda A_{n}\big(\blh,\blP\big)+B\big(\blh,\blP\big)
+2\zeta_n+R_n+\cR_n\Big\},
\end{eqnarray}
where we have used also $\sqrt{3}<2$.

$\mathbf{3^0}.$ Note  that in view of $\blP\in\mP(f)$, (\ref{eq3:proof-th:sup-norm-oracle}) and (\ref{eq9:proof-th:sup-norm-oracle})
\begin{eqnarray}
\label{eq10:proof-th:sup-norm-oracle}
\Big\|\widehat{f}_{\blh,\blP}-f\Big\|_\infty&=&\bigg\|\prod_{\blI\in\blP}
\widetilde{f}_{\blh_{\blI}}\big(x_{\blI}\big)-\prod_{\blI\in\blP}f_{\blI}(x_{\blI})\bigg\|_\infty
\nonumber\\
&\leq&
d\left[\max\big\{\bar{\mathbf{f}}_n,\rk^{2}_1\mathbf{f}\big\}\right]^{d-1}
\sup_{\blI\in\blP}\left\|\widetilde{f}_{\blh_{\blI}}\big(x_{\blI}\big)-f_{\blI}(x_{\blI})\right\|_{\blI,\infty}
\nonumber\\
&\leq& d\left[\max\big\{\bar{\mathbf{f}}_n,\rk^{2}_1\mathbf{f}\big\}\right]^{d-1}\Big[B(\blh,\blP)+\zeta(\blh,\blP)\Big]
\leq \bar{\mathfrak{f}}_n\Big[B(\blh,\blP)+\Lambda A_n(\blh,\blP)+\zeta_n  \Big].
\end{eqnarray}
Here we have also used that $\cP\equiv\cP\diamond\cP$.
We obtain from (\ref{eq1:proof-th:sup-norm-oracle}),   (\ref{eq99:proof-th:sup-norm-oracle}) and (\ref{eq10:proof-th:sup-norm-oracle})
\begin{eqnarray*}
\Big\|\widehat{f}_{\widehat{h},\widehat{\cP}}-f\Big\|_\infty\leq
\bar{\mathfrak{f}}_n\left[3B\big(\blh,\blP\big)+7\Lambda A_{n}(\blh,\blP)
+5\zeta_n+2R_n+2\cR_n\right],
\end{eqnarray*}
and, therefore, for any  $\blh\in\cH_n$,  $\blP\in\mP$ and $q\geq 1$
\begin{eqnarray}
\label{eq11:proof-th:sup-norm-oracle}
&&\;\left(\bE^{(n)}_f\big\|\widehat{f}_{\widehat{h}[\widehat{\cP}],\widehat{\cP}}-f\big\|_\infty\right)^{\frac{1}{q}}\leq
E_q\Big[3B\big(\blh,\blP\big)+7\Lambda A_{n}(\blh,\blP)\Big]+E_{2q}\Big[5y_{1,n}+2\Lambda\ma^* \big(y_{2,n}+y_{3,n}\big)\Big],
\end{eqnarray}
where we have put for $p\geq 1$
$$
E_p=\left[\bE^{(n)}_f\big(\bar{\mathfrak{f}}_n\big)^p\right]^{\frac{1}{p}},
y_{1,n}=\left[\bE^{(n)}_f\big(\zeta_n\big)^{2q}\right]^{\frac{1}{2q}},
y_{2,n}=\left[\bE^{(n)}_f\left[\bar{s}_n-\bar{\mathbf{f}}_n\right]_+^{2q}\right]^{\frac{1}{2q}},
y_{3,n}=\left[\bE^{(n)}_f\left[\bar{\mathbf{f}}_n-3\bar{s}_n\right]_+^{2q}\right]^{\frac{1}{2q}}.
$$
Taking into account that the right hand side of (\ref{eq11:proof-th:sup-norm-oracle}) is independent of the choice $\blh$ and $\blP$
and that the quantity $\bar{s}_n\leq 1\vee [\rk_1\mathbf{f}]$
 we get
\begin{eqnarray*}
\label{eq12:proof-th:sup-norm-oracle}
\left(\bE^{(n)}_f\big\|\widehat{f}_{\widehat{h}[\widehat{\cP}],\widehat{\cP}}-f\big\|_\infty\right)^{\frac{1}{q}}&\leq&
7\Lambda E_q\bigg(\inf_{\blh\in\cH_n}\inf_{\blP\in\mP(f)}\Big[B\big(\blh,\blP\big)+A_{n}(\blh,\blP)\Big]\bigg)
\\
&\qquad& +E_{2q}\Big[5y_{1,n}+2\Lambda\ma^* \big(y_{2,n}+y_{3,n}\big)\Big]
\\
&=& \mathbf{C_1}\big(q,d,\mathbf{K},\mathbf{f}\big) \mathfrak{R}(f) +E_{2q}\Big[5y_{1,n}+2\Lambda\ma^* \big(y_{2,n}+y_{3,n}\big)\Big].
\end{eqnarray*}
where we have put $\mathbf{C_1}\big(q,d,\mathbf{K},\mathbf{f}\big)=7\Lambda E_q\sqrt{1\vee [\rk_1\mathbf{f}]}$.

This inequality together with  bounds found in  Lemma \ref{lem:stoch-part} leads to the assertion of the theorem.

\epr

\subsection{Proof of Theorem \ref{th:lower-bound-in-supnorm}}

The proof of Theorem \ref{th:lower-bound-in-supnorm} is relatively standard and  based on the  general result established in
\cite{lepski-kerk2}, Proposition 7. For the convenience we formulate  this result not in  full generality but its  version reduced to the considered problem.
 Let $\big(\beta,p,\cP\big)\in (0,\infty)^{d}\times[1,\infty]^d\times\mP$ such that $\Upsilon\big(\beta,p,\cP\big)>0$ and  $\cL\in (0,\infty)^d$ be fixed.

\begin{lemma}
\label{lem:KLP-result}
Assume that there exist $f_0\in \mathbf{N}_{p,d}\big(\beta,\cL,\cP\big),$ $ \rho_n>0,\;n\in\bN^*$, and a finite  set $\mfJ_n$  such that for any sufficiently large $n\in\bN^*$   one can find
 $\left\{f^{(\mfj)},\;\mfj\in\mfJ_n\right\}\subset \mathbf{N}_{p,d}\big(\beta,\cL,\cP\big)$ satisfying
\begin{eqnarray}
\label{eq:ass1-klp-lemma}
  &&\|f^{(\mfj)}- f_0\|_\infty =\rho_n,\qquad\; \forall \mfj\in\mfJ_n;
\\*[2mm]
\label{eq:ass2-klp-lemma}
&& \limsup_{n\to\infty}\bE^{(n)}_{f_0}\Bigg[\frac{1}{|\mfJ_n|}\sum_{\mfj\in\mfJ_n}\frac{\rd \bP^{(n)}_{f^{(\mfj)}}}{\rd \bP^{(n)}_{f_0}}\Big(X^{(n)}\Big)-1\Bigg]^{2}=:\mathbf{C}<\infty.
\end{eqnarray}
Then for $r\geq 1$
$$
\liminf_{n\to\infty}\inf_{\tilde f} \sup_{f \in \mathbf{N}_{p,d}\big(\beta,\cL,\cP\big)} \rho^{-1}_n\left(\bE^{(n)}_f \big\|\tilde{f} - f\big\|^{r}_\infty\right)^{\frac{1}{r}}
\geq 2^{-1} \left[1-\sqrt{\mathbf{C}/(\mathbf{C}+4)}\right],
$$
where infimum is taken over all possible estimators.

\end{lemma}

\paragraph{Proof of the theorem}

Set
$
\mathcal{N}(x)=\prod_{i=1}^d\bigg(\left[2\pi\right]^{-1/2}\exp{-\big\{x^{2}_i/2\big\}}\Bigg)
$
and let $f_0(x)=\sigma^{-1}\cN(x/\sigma)$,
where $\sigma>0$ is chosen in such a way that $f_0$ belongs to the class $\overline{\bN}_{p,d}\big(\beta,\cL/2\big)$.
We remark that in order to obey the latter restriction it suffices to choose $\sigma$ satisfying
\begin{equation}
\label{eq0:proof-th:lower-bound-in-supnorm}
\sup_{\mfI\in\cI_d}\sup_{i\in\mfI}\sigma^{-\beta_i+|\mfI|/r_i}\big\|\cN_{\mfI}\big\|^{|\mfI|}_{r_i}\leq\inf_{i=\overline{1,d}}\cL_i.
\end{equation}
The product structure of  $f_0$ together with (\ref{eq0:proof-th:lower-bound-in-supnorm}) allows us to assert that $f_0\in \mathbf{N}_{p,d}\big(\beta,\cL/2,\cP\big)$ for any $\cP\in\mP$.
Let $\mfI^*\in\{1,\ldots,d\}$ be defined from the relation
$$
\Upsilon\big(\beta,p,\cP\big):=\inf_{\mfI\in\cP}\gamma_\mfI(\beta,p)=\gamma_{\mfI^*}(\beta,p),
$$
and for the notation convenience the elements of $\mfI^*$ will be denoted by $i_1,\ldots,i_m$ and $m=|I^*|$.

Let  $g:\bR\to\bR$ be  compactly supported on $(-1/2,1/2)$ function, satisfying $g\in\cap_{i\in\mfI^*}N_{p_i,1}(\beta_i,1/2)$, and such that
$
\int g=0.
$
Suppose also that $\big|g(0)\big|=\|g\|_\infty$.

\vskip0.1cm

Let $A_n\to 0$ and $\delta_{l,n}\to 0,\;l=\overline{1,m},$  $n\to \infty$, be  sequences whose choice will be done later and set
$\mfJ_n:=\big[1,\ldots, M_{1,n}\big]\times\cdots\times\big[1,\ldots,M_{m,n}\big]\subset\bN^{m}$, where $M_{l,n}=\big\lfloor\delta^{-1/2}_{l,n}\big\rfloor,\;l=\overline{1,m}.$

\vskip0.1cm

For any $\mfj=\big(j_1,\ldots j_m\big)\in\mfJ_n$ define
$
G_{\mfj}\big(x_{\mfI})=A_n\prod_{l=1}^m g\left(\delta^{-1}_{i,n}\Big[x_{i_l}-x^{(\mfj)}_{i_l}\Big]\right).
$
Here  for any $\mfj\in\mfJ_n$ we put $x^{(\mfj)}_{i_l}=j_l\delta_{l,n}$. The choice of $g$  implies
\begin{equation}
\label{eq1:proof-th:lower-bound-in-supnorm}
G_{\mfj}G_{\mfj^\prime}\equiv 0,\;\;\;\forall \mfj,\mfj^\prime\in\mfJ_n,\;\mfj\neq\mfj^\prime.
\end{equation}
Note also that the system of equations
\begin{equation}
\label{eq2:proof-th:lower-bound-in-supnorm}
A_n\delta^{-\beta_{i_k}}_{k,n}\bigg(\prod_{l=1}^m\delta_{l,n}\bigg)^{1/p_{i_k}}=\frac{\cL_{i_k}}{c_k},\;\; k=\overline{1,m},
\end{equation}
implies that $G_{\mfj}\in N_{p_{\mfI},d}\big(\beta_{\mfI},\cL_{\mfI}\big/2\big)$ for any $\mfj\in\mfJ_n$.
Here we have denoted $c_k=\big(\|g\|_{p_{i_k}}\big)^{m-1}$.

Introduce the family of functions
$
\left\{f^{(\mathbf{j})},\;\mathbf{j}\in\mfJ_n\right\}
$
as follows.
$$
f^{(\mathbf{j})}(x)=\prod_{i\notin \mfI^*}^d\bigg(\left[2\pi\sigma^2\right]^{-1/2}\exp{-\big\{x^{2}_i/2\sigma^2\big\}}\Bigg)
\bigg(\prod_{i\in\mfI^*}^d\left[2\pi\sigma^2\right]^{-1/2}\exp{-\big\{x^{2}_i/2\sigma^2\big\}}+G_{\mfj}\big(x_{\mfI})\Bigg).
$$
First we remark that  $A_n\to 0,\;n\to\infty$, implies that $f^{(\mathbf{j})}>0$ for all sufficiently large $n$. Next, the assumption
$\int g=0$ implies that $\int f^{(\mathbf{j})}=1$. Thus,  $f^{(\mathbf{j})}$ is a probability
density for any $\mfj\in\mfJ_n$ for all sufficiently large $n$.
At last the choice of $f_0$ together with (\ref{eq2:proof-th:lower-bound-in-supnorm}) allows us to assert that $f^{(\mathbf{j})}\in \mathbf{N}_{p,d}\big(\beta,\cL,\cP\big)$  for any $\mfj\in\mfJ_n$.

Thus, we conclude that Lemma \ref{lem:KLP-result} is applicable to the family
$
\left\{f^{(\mathbf{j})},\;\mathbf{j}\in\mfJ_n\right\}.
$
We remark also that
\begin{equation}
\label{eq3:proof-th:lower-bound-in-supnorm}
\big\|f^{(\mathbf{j})}-f_0\big\|_\infty= c^*_1A_n,\;\;\forall \mathbf{j}\in\mfJ_n,
\end{equation}
where we have put $c^*_1=|g(0)|^{m}\Big(2\pi\sigma^2\Big)^{(m-d)/2}$. Here we have also used that $\big|g(0)\big|=\|g\|_\infty$.
We conclude that the assumption (\ref{eq:ass1-klp-lemma}) is fulfilled  with $\rho_n=c^*_1A_n$.

Let us now proceed with the verification of the condition (\ref{eq:ass2-klp-lemma}) of Lemma \ref{lem:KLP-result}. Note first that
$$
\frac{\rd \bP^{(n)}_{f^{(\mfj)}}}{\rd \bP^{(n)}_{f_0}}\Big(X^{(n)}\Big)=\prod_{k=1}^n\frac{f^{(\mathbf{j})}(X_k)}{f_0(X_k)}
$$
and, therefore,
\begin{equation}
\label{eq4:proof-th:lower-bound-in-supnorm}
\Bigg[\frac{1}{|\mfJ_n|}\sum_{\mfj\in\mfJ_n}\frac{\rd \bP^{(n)}_{f^{(\mfj)}}}{\rd \bP^{(n)}_{f_0}}\Big(X^{(n)}\Big)\Bigg]^{2}=
\frac{1}{|\mfJ_n|^{2}}\Bigg\{\sum_{\mfj\in\mfJ_n}\prod_{k=1}^n\bigg[\frac{f^{(\mathbf{j})}(X_k)}{f_0(X_k)}\bigg]^{2}+
\sum_{\substack{\mfj,\mfj^{\prime}\in \mfJ_n:\\
\mfj\neq\mfj^{\prime}}}\prod_{k=1}^n\frac{f^{(\mathbf{j})}(X_k)f^{(\mathbf{j^\prime})}(X_k)}{f^{2}_0(X_k)}
\Bigg\}.
\end{equation}
Since  $X_k,\; k=\overline{1,n}$ are i.i.d. random vectors, we have for any $\mfj\neq\mfj^\prime$
\begin{equation*}
\bE^{(n)}_{f_0}\Bigg\{
\prod_{k=1}^n\frac{f^{(\mathbf{j})}(X_k)f^{(\mathbf{j^\prime})}(X_k)}{f^{2}_0(X_k)}
\Bigg\}=\Bigg\{
\int_{\bR^{|\mfI^*|}}\Bigg[1+\frac{G_{\mfj}\big(x_{\mfI^*}\big)}{f_{\mfI^*,0}\big(x_{\mfI^*}\big)}\Bigg]
\Bigg[1+\frac{G_{\mfj^\prime}\big(x_{\mfI^*})}{f_{\mfI^*,0}\big(x_{\mfI^*}\big)}\Bigg]f_{\mfI^*,0}\big(x_{\mfI^*}\big)\rd x_{\mfI^*}
\Bigg\}^{n}
=1.
\end{equation*}
To get the last equality we have used (\ref{eq1:proof-th:lower-bound-in-supnorm}) and the fact that $\int_{\bR^{|\mfI^*|}}G_{\mfj}\big(x_{\mfI^*})\rd x_{\mfI^*}=0$ since $\int g=0$.

The latter result together with (\ref{eq4:proof-th:lower-bound-in-supnorm}) yields
\begin{eqnarray}
\label{eq5:proof-th:lower-bound-in-supnorm}
\cE_n&:=&\bE^{(n)}_{f_0}\Bigg[\frac{1}{|\mfJ_n|}\sum_{\mfj\in\mfJ_n}\frac{\rd \bP^{(n)}_{f^{(\mfj)}}}{\rd \bP^{(n)}_{f_0}}\Big(X^{(n)}\Big)-1\Bigg]^{2}
\nonumber\\
&=&
\frac{1}{|\mfJ_n|^{2}}\sum_{\mfj\in\mfJ_n}\Bigg\{
\int_{\bR^{|\mfI^*|}}\Bigg[1+\frac{G_{\mfj}\big(x_{\mfI^*}\big)}{f_{\mfI^*,0}\big(x_{\mfI^*}\big)}\Bigg]^{2}
f_{\mfI^*,0}\big(x_{\mfI^*}\big)\rd x_{\mfI^*}
\Bigg\}^{n}-|\mfJ_n|^{-1}
\nonumber\\
&=&\frac{1}{|\mfJ_n|^{2}}\sum_{\mfj\in\mfJ_n}\Bigg\{
1+\int_{\bR^{m}}\bigg[\frac{G^{2}_{\mfj}(y)}{f_{\mfI^*,0}(y)}\bigg]
\rd y
\Bigg\}^{n}-|\mfJ_n|^{-1}.
\end{eqnarray}
Since, $G_{\mfj}(y)=0$ for any $y\notin \big[0,\sqrt{\delta_{1,n}}\big]\times\cdot\times\big[0,\sqrt{\delta_{m,n}}\big]=:\cY_n$ we have for all
$n$ large enough
$
\inf_{y\in\cY_n}f_{\mfI^*,0}(y)\geq 2^{-1}\Big(2\pi\sigma^2\Big)^{-m}.
$
It yields together with (\ref{eq4:proof-th:lower-bound-in-supnorm}), putting $c^*_2=2\Big(2\pi\sigma^2\Big)^{m}\|g\|^{2m}_2$,
$$
\cE_n\leq |\mfJ_n|^{-1}\bigg(
1+c^*_2A^{2}_n\prod_{l=1}^m\delta_{l,n}
\bigg)^{n}.
$$
If we choose $A_n$ and $\delta_{l,n},\;l=\overline{1,m}$ satisfying
\begin{eqnarray}
\label{eq7:proof-th:lower-bound-in-supnorm}
c^*_2nA^{2}_n\prod_{l=1}^m\delta_{l,n}\leq (1/4)\ln\Big(\prod_{l=1}^m\delta^{-1}_{l,n}\Big)\leq \ln{\big(|\mfJ_n|\big)},
\end{eqnarray}
for all $n\geq 1$ large enough, then $\cE_n\leq 1$ and, therefore, the condition (\ref{eq:ass2-klp-lemma}) is fulfilled with $\mathbf{C}=1$.

Thus, we have to choose $A_n$ and $\delta_{l,n},\;l=\overline{1,m}$ satisfying (\ref{eq2:proof-th:lower-bound-in-supnorm}) and (\ref{eq7:proof-th:lower-bound-in-supnorm}). Let $t>0$ be the number whose choice will be done later. Consider instead of (\ref{eq7:proof-th:lower-bound-in-supnorm}) the equation
\begin{eqnarray}
\label{eq11:proof-th:lower-bound-in-supnorm}
nA^{2}_n\prod_{l=1}^m\delta_{l,n}= t^{2}\ln(n).
\end{eqnarray}
and solve (\ref{eq2:proof-th:lower-bound-in-supnorm}) and (\ref{eq11:proof-th:lower-bound-in-supnorm}). Straightforward computations yield
$$
A_n=R(\e t)^{\frac{1-\sum_{l=1}^m \frac{1}{\beta_{i_l}p_{i_l}}}{1-\sum_{l=1}^m \left(\frac{1}{p_{i_l}}-\frac{1}{2}\right)\frac{1}{\beta_{i_l}}}},\qquad
\delta_{l,n}=A_n^{\frac{1}{\beta_{i_l}}-\frac{2}{\beta_{i_l}p_{i_l}}}\;\big(t\e)^{\frac{2}{\beta_{i_l}p_{i_l}}}
\big(c_l/\cL_l\big)^{\frac{1}{\beta_{i_l}}},
$$
where we have put $R=\left(\prod_{l=1}^m\big(c_l/\cL_l\big)^{\frac{1}{2\beta_{i_l}}}\right)^{\frac{1}{1-\sum_{l=1}^m \left(\frac{1}{p_{i_l}}-\frac{1}{2}\right)\frac{1}{\beta_{i_l}}}}$. Moreover we have in view of  (\ref{eq11:proof-th:lower-bound-in-supnorm})
$$
\left(\prod_{l=1}^m\delta_{l,n}\right)^{-1/2}=R(\e t)^{-a},\quad a=\frac{\sum_{l=1}^m \frac{1}{\beta_{i_l}}}{1-\sum_{l=1}^m \left(\frac{1}{p_{i_l}}-\frac{1}{2}\right)\frac{1}{\beta_{i_l}}}
$$
and, therefore,
$
(1/4)\ln\Big(\prod_{l=1}^m\delta^{-1}_{l,n}\Big)\asymp
(a/2)\ln(n),\;\; n\to\infty,
$
Hence, choosing $t$ as an arbitrary number satisfying
$
t^{2}< (2c^*_2)^{-1}a
$
we guarantee that (\ref{eq11:proof-th:lower-bound-in-supnorm}) implies
(\ref{eq7:proof-th:lower-bound-in-supnorm}) for all $n$ large enough.

Thus, we conclude that Lemma \ref{lem:KLP-result} is applicable with
$$
\rho_n=c^*_1A_n=c^*_1 R\left(\frac{t\ln(n)}{n}\right)^{\frac{1-\sum_{l=1}^m \frac{1}{\beta_{i_l}p_{i_l}}}{2\Big(1-\sum_{l=1}^m \Big[\frac{1}{p_{i_l}}-\frac{1}{2}\Big]\frac{1}{\beta_{i_l}}\Big)}}.
$$
It remains to note that the definition of $I^*$ implies that
$
\Upsilon\big(\beta,p,\cP\big)=\frac{1-\sum_{l=1}^m\frac{1}{\beta_{i_l}p_{i_l}}}{\sum_{l=1}^m\frac{1}{\beta_{i_l}}}.
$
We remark that
$$
\frac{\Upsilon\big(\beta,p,\cP\big)}{2\Upsilon\big(\beta,p,\cP\big)+1}=\frac{1-\sum_{l=1}^m \frac{1}{\beta_{i_l}p_{i_l}}}{2\left(1-\sum_{l=1}^m \left[\frac{1}{p_{i_l}}-\frac{1}{2}\right]\frac{1}{\beta_{i_l}}\right)}
$$
and the assertion of the theorem follows.
\epr

\subsection{Proof of Theorem \ref{th:adaptation-in-supnorm}}

The proof of the theorem is based on the application of Theorem \ref{th:sup-norm-oracle}  and on Lemma \ref{lem:bound-for-bais-supnorm} below
that allows us to bound from above the quantity $B(h,\cP)$.
The assertion of the lemma, whose proof is postponed to Appendix, is based on the embedding  theorem for anisotropic Nikolskii spaces.
For any function $g:\bR^s\to\bR$ and any $\eta\in (0,\infty)^s$ set
$$
\cB_{\eta,g}(z)=\int_{\bR^s}K_\eta(t-z)g(t)\rd t- g(z),\;\; z\in\bR^s.
$$

\begin{lemma}
\label{lem:bound-for-bais-supnorm} Let $\mathbf{K}$ satisfy Assumption \ref{ass:on-kernel} and (\ref{eq:ass-th:adaptation-in-supnorm}).
Let $(\alpha,r)\in (0,\mb]^{s}\times[1,\infty]^s$ be such that
$
\kappa=1-\sum_{l=1}^{s} (\alpha_l r_l)^{-1}>0
$
and let $Q\in (0,\infty)^s$.
Then  there exists $\mathrm{c}=\mathrm{c}\big(s,r,\mb\big)>0$   such that
$$
\sup_{g\in\bN_{r, s}(\alpha,Q)}\left\|\cB_{\eta,g}\right\|_\infty\leq \mathrm{c}\rk_1^{s}\sum_{i=1}^s Q_i\eta^{\blalpha_i}_i,\;\;\forall \eta\in (0,\infty)^s.
$$
Here $\blalpha=(\blalpha_1,\ldots\blalpha_s\;)$, $\blalpha_i=\kappa\alpha_i\kappa_i^{-1}$ and
$\kappa_i=1-\sum_{l=1}^{s}\big(r^{-1}_l-r^{-1}_i\big)\alpha_l^{-1}$.

\end{lemma}

\paragraph{Proof of Theorem \ref{th:adaptation-in-supnorm}} Let $\big(\beta,p,\cP\big)\in (0,\mb]^{d}\times[1,\infty]^d\times\mP$ such that $\Upsilon\big(\beta,p,\cP\big)>0$ and  $\cL\in (0,\infty)^d$  be fixed. For any $\mfI\in\cI_d$  and any $\mfi\in\mfI$ define
$$
\blb_\mfi(\mfI)=\tau(\mfI)\beta_\mfi\tau_\mfi^{-1}(\mfI),\quad\tau(\mfI)=1-\sum_{l\in\mfI} (\beta_l p_l)^{-1},\quad
\tau_\mfi(\mfI)=1-\sum_{l\in\mfI}\big(p^{-1}_l-p^{-1}_\mfi\big)\beta_l^{-1},
$$
and  remark that the condition $\Upsilon\big(\beta,p,\cP\big)>0$ implies that $\tau(\mfI)>0$ for any $\mfI\in\cI_d$.

\smallskip

Let us first prove the following simple fact. Denote  $\cC_\mfi(\mfI)=\{\mathbf{J}\subseteq\mfI:\;\mfi\in \mathbf{J}\},\;\mfi\in\mfI$. Then \begin{equation}
\label{eq1:proof-th:adaptation-in-supnorm}
\blb_\mfi(\mfI)=\inf_{\mathbf{J}\in\cC_\mfi(\mfI)} \blb_\mfi(\mfJ),\;\;\forall \mfi\in\mfI.
\end{equation}
Indeed, we remark that
$
\tau_\mfi(\mfJ)=1-\sum_{l\in\mfJ}\big(p^{-1}_l-p^{-1}_\mfi\big)\beta_l^{-1}=\tau(\mfJ)+p^{-1}_\mfi\sum_{l\in\mfJ}\beta_l^{-1}
$
and, therefore,
$$
\blb_\mfi(\mfJ)=\frac{\beta_\mfi\tau(\mfJ)}{\tau(\mfJ)+p^{-1}_\mfi \beta^{-1}(\mfJ)},\qquad \beta^{-1}(\mfJ)=\sum_{l\in\mfJ}\beta_l^{-1}.
$$
We obviously have $\tau(\mfJ)\geq\tau(\mfI)$ and $\beta^{-1}(\mfJ)\leq \beta^{-1}(\mfI)$ for any $\mfJ\subseteq\mfI$. It remains to note that
 $x\mapsto x/(x+a)$ is increasing on $\bR_+$ for any $a>0$ and (\ref{eq1:proof-th:adaptation-in-supnorm})
follows.

\smallskip

Let $\cP^{\prime}\in\mP$ be an arbitrary partition. Since $f\in\overline{\bN}_{p,d}\big(\beta,\cL\big)$ we have
$
 f_{\mfJ}\in \bN_{p_\mfJ,|\mfJ|}\big(\beta_\mfJ,\cL_\mfJ\big)
$
 and, therefore, in view of Lemma \ref{lem:bound-for-bais-supnorm} we have for any $h\in (0,1]^d$ and $\mfJ\in\cP\diamond\cP^\prime$
$$
b_{h_{\mfJ}}\leq \mathrm{c}\big(|\mfJ|,p_{\mfJ},\mb\big)\rk_1^{|\mfJ|}\sum_{\mfi\in\mfJ} \cL_\mfi h^{\blb_\mfi(\mfJ)}_\mfi
\leq \mathrm{c}_1\sum_{\mfi\in\mfI} \cL_\mfi h^{\blb_\mfi(\mfI)}_\mfi.
$$
To get the last inequality we use (\ref{eq1:proof-th:adaptation-in-supnorm}),  $h\in (0,1]^d$  and we have put
$\mathrm{c}_1=\rk_1^{d}\sup_{\mfJ\in\cI_d}\mathrm{c}\big(|\mfJ|,p_{\mfJ},\mb\big)\rk_1^{|\mfJ|}$.

Noting that the right hand side of the latter inequality is independent  on $\mfJ$ we obtain
$$
B\big(h,\cP\big)\leq  \mathrm{c}_1\sup_{\mfI\in\cP} \sum_{\mfi\in\mfI} \cL_\mfi h^{\blb_\mfi(\mfI)}_\mfi,\quad h\in (0,1]^d.
$$
It remains to choose multi-bandwidth  $h$.  To do it it suffices to solve  for any $\mfI\in\cP$ the following system of equations.
$$
\cL_\mfj h^{\blb_\mfj(\mfI)}_\mfj=\cL_\mfi h^{\blb_\mfi(\mfI)}_\mfi=\sqrt{\frac{\ln(n)}{nV_{h_\mfI}}}, \;\;\mfi,\mfj\in\mfI.
$$
The solution is given by
$$
h_\mfi=\cL^{-\frac{1}{\blb_\mfi(\mfI)}}\left(\frac{\mL(\mfI)\ln(n)}{n}\right)^{\frac{\gamma_{\mfI}(\beta,p)}{2+\gamma_{\mfI}(\beta,p)}},\quad
\mL(\mfI)=\prod_{\mfi\in\mfI}\cL_{\mfi}^{\frac{1}{\beta_\mfi(\mfI)}}.
$$
Here we have also used that $1/\gamma_{\mfI}(\beta,p)=\sum_{\mfi\in\mfI}1/\blb_\mfi(\mfI)$.
The assertion of the theorem follows now from Theorem \ref{th:sup-norm-oracle}.
\epr

\section{Appendix}
\label{sec:appendix}

\subsection{Proof of Proposition \ref{prop:exp-ineq-for-kernel-process}}

$1^0.\;$ Note  that  $\mathbf{M}(z)=\mathbf{M}(|z|)$ since  $\mathbf{M}$ is symmetric that implies
$$
\chi_r(y)=n^{-1}\sum_{j=1}^n\left[M_r\big(\vec{\rho}\left(Y_j,y\right)\big)-
\bE^{(n)}_g\Big\{M_r\big(\vec{\rho}(Y_j,y)\big)\Big\}\right],
$$
where $\vec{\rho}:\bR^s\times\bR^s\to\bR^s$ is given by $\vec{\rho}\big(z,z^\prime\big)=\big(|z_1-z_1^\prime|,\ldots,|z_s-z_s^\prime|\big)$.

We conclude that considered family of random fields obeys the structural assumption introduced in Section 4.4. of \cite{Lep2012}, with
$d=s$, $\bX_1^d=\bar{\bX}_1^d=\bR^s$ and $\rho_l:\bR\times\bR\to\bR$ is given by $|z-z^\prime|$ for any $l=\overline{1,s}$.
It implies in particular that $\bR^s$ is equipped with the metric $\varrho_s$ generated by the supremum norm, i.e.
$\varrho_s=\max_{l=\overline{1,s}}\rho_l$.
We remark also that in our case $K(u)=\prod_{l=1}^s\mathbf{M}(u_l),\;u\in\bR^s,\;$ $g\equiv 1$ and
$\gamma_l=1, l=\overline{1,s}$.

To get the assertion of Proposition \ref{prop:exp-ineq-for-kernel-process} we will apply Theorem 9 in \cite{Lep2012}
on $\cR_n(s):=[1/n, 1]^s$. Note that obviously  $\widetilde{\cR}_n\subseteq\cR_n(s)$.
Thus,  we have  to check the assumptions of the latter theorem and to match the notations used in the present paper and in \cite{Lep2012}.

\smallskip

First we note that since $\mathbf{M}$ satisfies Assumption \ref{ass:on-kernel}  Assumption 9 ($\mfi$) is obviously fulfilled with
$L_1=(3s/2)(\mathrm{m}_\infty)^{s-1}L$. Moreover Assumption 9 ($\mfi\mfi$) holds because $g\equiv 1$.

Thus, Assumption 9 is checked.

\smallskip

Consider the collection of closed cubs
$
\bB_{\frac{1}{2}}(\mathbf{j})=\left\{z\in\bR^s:\;\; \varrho_s(z,\mathbf{j})\leq 1\right\},\;\mathbf{j}\in\bZ^{s},
$
and let $\mathfrak{E}_{\mathbf{j}}(\delta),\;\delta>0$ denote the metric entropy of $\bB_{\frac{1}{2}}(\mathbf{j})$ measured in the metric
$\varrho_s$.

Obviously $\left\{\bB_{\frac{1}{2}}(\mathbf{j}),\;\mathbf{j}\in\bZ^{d}\right\}$ is a countable cover of $\bR^s$ and each member of this collection
is totally bounded (even compact) subset of $\bR^s$.
It is easily seen that
$$
\text{card}\left(\left\{\mathbf{k}\in\bZ^{s}:\; \bB_{\frac{1}{2}}(\mathbf{j})\cap\bB_{\frac{1}{2}}(\mathbf{k})\neq\emptyset\right\}\right)\leq 3^s,\;\;\forall \mathbf{j}\in\bZ^{s}.
$$
Using the terminology of \cite{Lep2012} we can say that $\left\{\bB_{\frac{1}{2}}(\mathbf{j}),\;\mathbf{j}\in\bZ^{d}\right\}$ is $3^s$-totally bounded cover of $\bR^s$.
Moreover,
$
\mathfrak{E}_{\mathbf{j}}(\delta)=s\big[\ln(1/\delta)\big]_+
$
for any $\delta>0$ and any $\mathbf{j}\in\bZ^s$.
All saying above allows us to assert that Assumption 7 ($\mfi$) is fulfilled with $\mathbf{I}=\bZ^s,$
$\mathrm{X}_{\mathbf{j}}=\bB_{\frac{1}{2}}(\mathbf{j})$,  $N=1.5s$ and $R=1$. It remains to note that Assumption 7 ($\mfi\mfi$) is automatically fulfilled  in our case since
$g\equiv 1$.

\smallskip

Also we note that for any $\mathbf{j},\mathbf{k}\in\bZ^s$ satisfying $\bB_{\frac{1}{2}}(\mathbf{j})\cap\bB_{\frac{1}{2}}(\mathbf{k})=\emptyset$
one has
$$
\inf_{x\in\bB_{\frac{1}{2}}(\mathbf{j})}\inf_{y\in\bB_{\frac{1}{2}}(\mathbf{k})}\varrho_s(x,y)\geq 1
$$
and, therefore, Assumption 11 is checked with $\mathfrak{t}=1$. At last we have for any $n\geq 1$
$$
\sup_{r\in\cR_n(s)}\sup_{u\notin (0,1]^s}\bigg|\prod_{l=1}^s\mathbf{M}(u_l/r_l)\bigg|=0,
$$
since $\text{supp}(\mathbf{M})\subseteq [-1/2,1/2]$. Hence, the condition (4.24) of Theorem 9 is fulfilled  as well that completes the verification of the assumptions of the theorem.

\smallskip

$2^0.\;$ Let us match the notations. First, in our case   $\mathbf{n_1}=\mathbf{n_2}=n$. Since $Y_j,\;j\geq 1,$ are identically distributed the quantity
denoted $F_{\mathbf{n_2}}\big(r,\bar{x}^{(d)}\big)$ is given now by
$G(r,y)=\int_{\bR^s}|M_r(x-y)|\mathbf{g}(x)\rd x$ and, therefore, is independent on $n$.
Here we have taken into account that $\bar{x}^{(d)}\in\bX^d=\bR^s$.

It is easily seen that
\begin{equation}
\label{eq1:proof-prop-exp-ineq-for-kernel-process}
G_n:=\sup_{r\in [1/n,1]^s}\|G(r,\cdot)\|_\infty\leq \min\Big[\mathrm{m}^{s}_1\|\mathbf{g}\|_\infty,  \mathrm{m}^{s}_\infty n^{s}\Big].
\end{equation}
 It yields, in particular, that $F_{\mathbf{n_2}}=G_{n}\leq \mathrm{m}^{s}_1\|\mathbf{g}\|_\infty $ for any $n\geq 1$.

Choosing in Theorem 9  $q=p,\;v=2p+2$, $z=1$ and remembering that $\bar{x}^{(d)}=y$, we have
$$
\widehat{\cU}^{(v,z,p)}\big(n,r,\bar{x}^{(d)}\big)\leq \gamma_p\big(s,\mathrm{m}_\infty\big)\sqrt{\frac{\bar{G}(r)\ln(n)}{nV_r}},
$$
for any $\bar{x}^{(d)}=y\in\bR^s$ and any $r\in\widetilde{\cR}_n(s)\subseteq\cR_n(s)$.
To get this assertion we have used that
$G_n\leq (\mathrm{m}_\infty n)^s$ in view of (\ref{eq1:proof-prop-exp-ineq-for-kernel-process}).

At last, taking into account that the right hand side of the latter inequality is independent on $y$, we deduce from Theorem 9 that
for any $p\geq 1$
$$
\bE^{(n)}_\mathbf{g}\Bigg\{\sup_{r\in\widetilde{\cR}_n(s)}\Big[\big\|\chi_r\big\|_\infty-\gamma_p\big(s,\mathrm{m}_\infty\big)
\sqrt{\frac{\bar{G}(r)\ln(n)}{nV_r}}\Bigg\}^p_+\leq c_1(p,s)\big[1\vee\mathrm{m}^{s}_1\|\mathbf{g}\|_\infty\big]^{\frac{p}{2}}n^{-\frac{p}{2}}+c_2(p,s)n^{-p},
$$
where  $c_1(p,s)=2^{7p/2+5}3^{p+5s+4}\Gamma(p+1)\pi^{p}\big(s,\mathrm{m}_\infty\big)$ and $c_2(p,s)=2^{p+1}3^{5s}$. Here we have also used
that $G_n\leq \mathrm{m}^{s}_1\|\mathbf{g}\|_\infty$ in view of (\ref{eq1:proof-prop-exp-ineq-for-kernel-process}) that implies
$\widehat{F}_{\mathbf{n_2}}\leq  1\vee\mathrm{m}^{s}_1\|\mathbf{g}\|_\infty$.
\epr

\subsection{Proof of Proposition \ref{prop2:exp-ineq-for-kernel-process}}

First, noting that $\gamma_p\big(s,\mathrm{m}_\infty\big)\sqrt{\ma}=1/2$
we obtain from Proposition \ref{prop:exp-ineq-for-kernel-process} that
\begin{equation}
\label{eq1:cor-prop-exp-ineq-for-kernel-process}
\bE^{(n)}_\mathbf{g}\bigg\{\sup_{r\in\widetilde{\cR}_n^{(a)}(s)}\bigg(\big\|\chi_r\big\|_\infty-
\frac{1}{2}\sqrt{\bar{G}(r)}\bigg)\bigg\}^p_+\leq c_n,
\end{equation}
where we have put for brevity  $c_n=c_1(p,s)\big[1\vee\mathrm{m}^{s}_1\|\mathbf{g}\|_\infty\big]^{\frac{p}{2}}n^{-\frac{p}{2}}+c_2(p,s)n^{-p}$.
Next, putting
$
\bar{\chi}_r(y)=\Upsilon_r(y)-\bE^n_\mathbf{g}\Upsilon_r(y)
$
we have in view if (\ref{eq1:cor-prop-exp-ineq-for-kernel-process})
\begin{equation}
\label{eq2:cor-prop-exp-ineq-for-kernel-process}
\bE^{(n)}_\mathbf{g}\bigg\{\sup_{r\in\widetilde{\cR}_n^{(a)}(s)}\bigg(\big\|\bar{\chi}_r\big\|_\infty-
\frac{1}{2}\sqrt{\bar{G}(r)}\bigg)\bigg\}^p_+\leq c_n.
\end{equation}
To get the latter result we remarked that if $\mathbf{M}$ satisfies Assumption \ref{ass:on-kernel} then $|\mathbf{M}|$ satisfies it as well and, therefore, Proposition \ref{prop:exp-ineq-for-kernel-process} is applicable to the process $\bar{\chi}_r(\cdot)$. It remains to note
that the function $\bar{G}(\cdot)$ is the same for both processes $\chi_r(\cdot)$ and $\bar{\chi}_r(\cdot)$.
We also note that
$$
G(r)=\displaystyle{\sup_{y\in\bR^s}}\left\{\bE^{(n)}_\mathbf{g}\Upsilon_r(y)\right\}
$$
and, therefore, for any $r\in (0,1]^s$ one has
\begin{equation}
\label{eq22:cor-prop-exp-ineq-for-kernel-process}
\bar{G}(r)=1\vee\left\|\bE^{(n)}_\mathbf{g}\Upsilon_r\right\|_\infty\leq 1\vee\left\|\Upsilon_r\right\|_\infty+\left\|\bar{\chi}_r\right\|_\infty,
\end{equation}
where we have used the obvious inequality
$
\big|||x||\vee||z||-||y||\vee||z||\big|\leq ||x-y||
$
being true for any  normed vector space.

Hence, putting $\zeta_n(\ma)=\sup_{r\in\cR_n^{(\ma)}(s)}\left[\left\|\bar{\chi}_r\right\|_\infty-\frac{1}{2}\sqrt{\bar{G}(r)}\right]_+$  we obtain for any
$r\in\cR_n^{(\ma)}(s)$
$$
\bar{G}(r)\leq \frac{1}{2}\sqrt{\bar{G}(r)}+1\vee\left\|\Upsilon_r\right\|_\infty+\zeta_n(\ma).
$$
It yields
$
\Big[\bar{G}(r)-2\left(1\vee\left\|\Upsilon_r\right\|_\infty\right)\Big]_+\leq 2\zeta_n(\ma)
$
and we have in view of  (\ref{eq2:cor-prop-exp-ineq-for-kernel-process})
$$
\bE^{(n)}_\mathbf{g}\bigg\{\sup_{r\in\cR_n^{(\ma)}(s)}\left[\bar{G}(r)-
2\left(1\vee\left\|\Upsilon_r\right\|_\infty\right)\right]\bigg\}^p_+\leq 2^{p}c_n.
$$
Similarly to (\ref{eq22:cor-prop-exp-ineq-for-kernel-process}) we have
\begin{equation*}
\label{eq3:cor-prop-exp-ineq-for-kernel-process}
1\vee\left\|\Upsilon_r\right\|_\infty\leq \bar{G}(r)+\left\|\bar{\chi}_r\right\|_\infty
\leq (3/2)\bar{G}(r)+\zeta_n(\ma)
\end{equation*}
and, therefore
$
\left[1\vee\left\|\Upsilon_r\right\|_\infty-(3/2)\bar{G}(r)\right]_+\leq\zeta_n(\ma).
$
Thus, we get from (\ref{eq2:cor-prop-exp-ineq-for-kernel-process})
$$
\bE^{(n)}_\mathbf{g}\bigg\{\sup_{r\in\cR_n^{(\ma)}(s)}\left[1\vee\left\|\Upsilon_r\right\|_\infty-(3/2)\bar{G}(r)\right]\bigg\}^p_+\leq c_n.
$$
\epr

\subsection{Proof of Lemma \ref{lem:bais-supnorm-result}} We have in view of Fubini theorem for any $x_\mfI\in\bR^{\mfI}$
\begin{eqnarray*}
s^*_{h_\mfI,\eta_\mfI}\big(x_\mfI\big)&=&\int_{\bR^{|\mfI|}}\left[K_{h_{\mfI}}\star K_{\eta_{\mfI}}\right]\big(t_\mfI-x_\mfI\big)f_{\mathbf{I}}
\big(t_{\mathbf{I}}\big)\rd t_\mfI=\int_{\bR^{|\mfI|}}\left[\int_{\bR^{|\mfI|}}K_{\eta_{\mfI}}\big(y_\mfI\big) K_{h_{\mfI}}\big(t_\mfI-x_\mfI-y_\mfI\big)\rd y_\mfI\right]f_{\mathbf{I}}
\big(t_{\mathbf{I}}\big)\rd t_\mfI
\\
&=&\int_{\bR^{|\mfI|}}K_{\eta_{\mfI}}\big(z_\mfI-x_\mfI\big)\left[\int_{\bR^{|\mfI|}} K_{h_{\mfI}}\big(t_\mfI-z_\mfI\big)
f_{\mathbf{I}}\big(t_{\mathbf{I}}\big)\rd t_\mfI\right]\rd y_\mfI
\\
&=& s_{h_\mfI}\big(x_\mfI\big)+\int_{\bR^{|\mfI|}}K_{\eta_{\mfI}}\big(z_\mfI-x_\mfI\big)\left[\int_{\bR^{|\mfI|}} K_{h_{\mfI}}\big(t_\mfI-z_\mfI\big)
\left\{f_{\mathbf{I}}\big(t_{\mathbf{I}}\big)-f_{\mathbf{I}}\big(z_{\mathbf{I}}\big)\right\}\rd t_\mfI\right]\rd z_\mfI.
\end{eqnarray*}
Therefore,
$$
\left\|s^*_{h_\mfI,\eta_\mfI}-s_{\eta_\mfI}\right\|_{\mfI,\infty}
\leq b_{h_{\mfI}} \int_{\bR^{|\mfI|}}\Big|K_{\eta_{\mfI}}\big(y_\mfI\big)\Big|\rd y_\mfI
\leq \mathrm{k}^{d}_1b_{h_{\mfI}}.
$$
\epr

\subsection{Proof of Lemma  \ref{lem:stoch-part}} The proof of the lemma is completely  based on application of Propositions \ref{prop:exp-ineq-for-kernel-process}--\ref{prop2:exp-ineq-for-kernel-process} and Corollary \ref{cor:prop:exp-ineq-for-kernel-process}.
\paragraph{Proof of $(\mfi)$} 
Remind that $\zeta(h,\cP)
=\displaystyle{\sup_{\mfI\in\cP}}\left\|\xi_{h_\mfI}\right\|_{\mfI,\infty}$ and
$$
\zeta_n=\displaystyle{\sup_{\eta\in\cH_n}\sup_{\cP\in\mP}}\bigg[\zeta\big(\eta,\cP\big)-\Lambda A_n(\eta,\cP)\bigg]_+.
$$
Then, we  have
\begin{equation}
\label{eq1:proof-lem:stoch-part}
\left[\bE^{(n)}_f\big(\zeta_n\big)^{2q}\right]^{\frac{1}{2q}}=\sum_{\cP\in\mP}\sum_{\mfI\in\cP}\Bigg(\bE^{(n)}_f\bigg\{
\sup_{\eta_{\mfI}\in\cH^{(\ma_{\mfi})}_n(|\mfI|)}\bigg[\left\|\xi_{\eta_\mfI}\right\|_{\mfI,\infty}-\gamma_{2q}\big(|\mfI|,\mathrm{k}_\infty\big) \sqrt{\frac{\bar{s}_n\ln(n)}{n V_{\eta_\mfI}}}\bigg]\bigg\}^{2q}_+\Bigg)^{\frac{1}{2q}},
\end{equation}
where we have put $\cH^{(\ma_{\mfi})}_n(|\mfI|)=\left\{\eta_{\mfI}\in (0,1]^{|\mfI|}:\;\; nV_{\eta_{\mfI}}\geq \ma_{\mfI}^{-1}\ln(n)\right\}$
and $\ma_{\mfI}=\left[2\gamma_{2q}\big(\mfI,\mathrm{k}_\infty\big)\right]^{-2}$.

\noindent To get the latter result we have used first that
$
 A_n(\eta,\cP)=\sup_{\mfI\in\cP}\sqrt{\frac{\bar{s}_n\ln(n)}{n V_{\eta_\mfI}}}
$
and the trivial inequality $\big[\sup_i x_i-\sup_iy_i\big]_+\leq \sup_i[x_i-y_i]_+$. Next we have used that
$
\Lambda=\sup_{\cP\in\mP}\sup_{\mfI\in\cP}\gamma_{2q}\big(|\mfI|,\mathrm{k}_\infty\big).
$
At last we have used that $\eta\in\cH_n$ implies $\eta_\mfI\in\cH^{(\ma_{\mfi})}_n(|\mfI|)$ for any $\mfI\in\cI_d$ in view of the definition of
$\ma^*$.

\smallskip

Note that for any  for any $\mfI\in\cI_d$ and any $\eta_\mfI\in (0,1]^{|\mfI|}$
$$
\bar{s}\geq 1\vee \bigg\|\int_{\bR^{\mfI}}\big|K_{\eta_{\mfI}}\big(t_\mfI-\cdot\big)\big|f_{\mathbf{I}}
\big(t_{\mathbf{I}}\big)\rd t_\mfI\bigg\|_{\mfI,\infty}=:\bar{F}_{\mfI}(\eta).
$$
We conclude that  Proposition \ref{prop:exp-ineq-for-kernel-process} is applicable with $\chi_r=\xi_{\eta_\mfI}$, $\mathbf{M}=\mathbf{K}$, $p=2q$, $s=|\mfI|$, $\ma=\ma_{\mfi}$, $\bar{G}=\bar{F}_\mfI$ and the assertion $(\mfi)$ follows with
$$
\mathbf{c_1}\big(2q,d,\mathbf{K},\mathbf{f}\big)=\sum_{\cP\in\mP}\sum_{\mfI\in\cP}
\left[c_1\big(2q,|\mfI|\big)\big[1\vee\mathrm{k}^{|\mfI|}_1\mathbf{f}\big]^{q}+c_2\big(2q,|\mfI|\big)\right].
$$

\paragraph{Proof of $(\mfi\mfi)$}
Put for any $ h\in\cH_n$ and $\mfI\in\cI_d$
$$
s_{\mfI}\big(h_{\mfi}\big)=\bigg\|\int_{\bR^{\mfI}}\big|K_{h_{\mfI}}\big(t_\mfI-\cdot\big)\big|f_{\mathbf{I}}
\big(t_{\mathbf{I}}\big)\rd t_\mfI\bigg\|_{\mfI,\infty}\;\quad
f_{\mfI,n}\big(h_{\mfi}\big)=\Big\|n^{-1}\sum_{i=1}^n\big|K_{h_{\mathbf{I}}}\left(X_{\mathbf{I},i}-\cdot\right)\big|\Big\|_{\mfI,\infty}.
 $$
We have similarly to (\ref{eq1:proof-lem:stoch-part})
$
\big[\bar{s}_n-\bar{\mathbf{f}}_n]_+\leq \displaystyle{ \sup_{\mfI\in\cI_d}\sup_{h_{\mfI}\in\cH^{(\ma_{\mfi})}_n(|\mfI|)}}\big[s_{\mfI}\big(h_{\mfI}\big)-2\mathbf{f}_{\mfI,n}\big(h_{\mfI}\big)]_+
$
and hence
\begin{equation*}
\label{eq2:proof-lem:stoch-part}
\left[\bE^{(n)}_f\left[\bar{s}_n-\bar{\mathbf{f}}_n\right]_+^{2q}\right]^{\frac{1}{2q}}\leq \sum_{\mfI\in\cI_d}\Bigg(\bE^{(n)}_f\bigg\{
\sup_{h_{\mfI}\in\cH^{(\ma_{\mfi})}_n(|\mfI|)}\Big[s_{\mfI}\big(h_{\mfI}\big)-2\mathbf{f}_{\mfI,n}\big(h_{\mfI}\big)\Big]\bigg\}^{2q}_+\Bigg)^{\frac{1}{2q}},
\end{equation*}
The assertion $(\mfi\mfi)$ follows now from the second statement of Proposition \ref{prop2:exp-ineq-for-kernel-process} with
$$
\mathbf{c_2}\big(2q,d,\mathbf{K},\mathbf{f}\big)=\sum_{\mfI\in\cI_d}
\left[c^\prime_1\big(2q,|\mfI|\big)\big[1\vee\mathrm{k}^{|\mfI|}_1\mathbf{f}\big]^{q}+c^\prime_2\big(2q,|\mfI|\big)\right].
$$
\paragraph{Proof of $(\mfi\mfi\mfi)$} We have
$
\big[\bar{\mathbf{f}}_n-3\bar{s}_n]_+\leq 2\displaystyle{ \sup_{\mfI\in\cI_d}\sup_{h_{\mfI}\in\cH^{(\ma_{\mfi})}_n(|\mfI|)}}\big[\mathbf{f}_{\mfI,n}\big(h_{\mfI}\big)-(3/2)s_{\mfI}\big(h_{\mfI}\big)]_+
$
and hence
\begin{equation*}
\label{eq22:proof-lem:stoch-part}
\left[\bE^{(n)}_f\left[\bar{\mathbf{f}}_n-3\bar{s}_n\right]_+^{2q}\right]^{\frac{1}{2q}}\leq 2\sum_{\mfI\in\cI_d}\Bigg(\bE^{(n)}_f\bigg\{
\sup_{h_{\mfI}\in\cH^{(\ma_{\mfi})}_n(|\mfI|)}\Big[\mathbf{f}_{\mfI,n}\big(h_{\mfI}\big)-(3/2)s_{\mfI}\big(h_{\mfI}\big)\Big]\bigg\}^{2q}_+\Bigg)^{\frac{1}{2q}},
\end{equation*}
The assertion $(\mfi\mfi\mfi)$ follows now from the first assertion of Proposition \ref{prop2:exp-ineq-for-kernel-process} with
$$
\mathbf{c_3}\big(2q,d,\mathbf{K},\mathbf{f}\big)=2\sum_{\mfI\in\cI_d}
\left[c_1\big(2q,|\mfI|\big)\big[1\vee\mathrm{k}^{|\mfI|}_1\mathbf{f}\big]^{q}+c_2\big(2q,|\mfI|\big)\right].
$$
\paragraph{Proof of $(\mfi\mathbf{v})$} Note that
\begin{eqnarray}
\label{eq3:proof-lem:stoch-part}
\bar{\mathfrak{f}}_n&:=&2d^{2}\rk_1^d\big(\bar{\mathbf{f}}_n\big)^{\left\lfloor d^{2}/4\right\rfloor+1}\left[\max\big\{\bar{\mathbf{f}}_n,\rk^{2}_1\mathbf{f}\big\}\right]^{d-1}
\nonumber\\
&\leq&
\beta \left[\big(\mathbf{f}_n\big)^{\left\lfloor d^{2}/4\right\rfloor+d}+\big(1+\rk^2_1\mathbf{f}\big)^{d-1}\big(\mathbf{f}_n\big)^{\left\lfloor d^{2}/4\right\rfloor+1}+\big(\mathbf{f}_n\big)^{d-1}+\big(1+\rk^2_1\mathbf{f}\big)^{d-1}\right],
\end{eqnarray}
where we have used $\rk_1\geq 1$ and   put $\beta=2d^{2}\rk_1^d 2^{\left\lfloor d^{2}/4\right\rfloor+d}$.
Thus, to get the assertion $(\mfi\mathbf{v})$ it suffices to bound from above $\bE_f\big(\mathbf{f}_n\big)^p,\;p\geq 1 $. We obviously have
\begin{eqnarray*}
\label{eq4:proof-lem:stoch-part}
\mathbf{f}_n\leq\sum_{\mfI\in\cI_d}\sup_{h_{\mfI}\in\cH^{(\ma_{\mfi})}_n(|\mfI|)}
\Big\|n^{-1}\sum_{i=1}^n\big|K_{h_{\mathbf{I}}}\left(X_{\mathbf{I},i}-\cdot\right)\big|\Big\|_{\mfI,\infty},
\end{eqnarray*}
and using Corollary \ref{cor:prop:exp-ineq-for-kernel-process} we get for $p\geq 1$
\begin{eqnarray}
\label{eq5:proof-lem:stoch-part}
\left[\bE^{(n)}_f\big(\mathbf{f}_n\big)^p\right]^{\frac{1}{p}}\leq\sum_{\mfI\in\cI_d}\big[1\vee\mathrm{k}^{|\mfI|}_1\mathbf{f}\big]^{\frac{1}{2}}
\left[\gamma_p\big(|\mfI|,\mathrm{k}_\infty\big)+\big\{c_1(p,|\mfI|)+c_2(p,|\mfI|)\big\}^{\frac{1}{p}}\right].
\end{eqnarray}
 The assertion $(\mfi\mathbf{v})$ follows now from (\ref{eq3:proof-lem:stoch-part}) and (\ref{eq5:proof-lem:stoch-part}).
 \epr

\subsection{Proof of Lemma \ref{lem:bound-for-bais-supnorm}}

The proof of the lemma is based on the  embedding  theorem for anisotropic Nikolskii classes which we formulate below.

Let $(\alpha,r)\in (0,\infty)^{s}\times[1,\infty]^s$ be such that
$
\kappa=1-\sum_{l=1}^{s} (\alpha_l r_l)^{-1}>0
$
and let $Q\in (0,\infty)^s$.
Then  there exists $\mathbf{c}>0$ completely determined by $\alpha,r$ and $s$ such that
\begin{equation}
\label{eq:inclusion-aniso-Nikolskii}
\bN_{r,s}\big(\alpha,Q\big)\subseteq \bN_{\infty,s}\big(\blalpha,\mathbf{c}Q\big),
\end{equation}
where $\blalpha=(\blalpha_1,\ldots\blalpha_s\;)$, $\blalpha_j=\kappa\alpha_j\kappa_j^{-1}$ and
$\kappa_j=1-\sum_{l=1}^{s}\big(r^{-1}_l-r^{-1}_j\big)\alpha_l^{-1}$.

The inclusion (\ref{eq:inclusion-aniso-Nikolskii}) is a particular case of Theorem 6.9 in \cite{nikol}, with $p^\prime=\infty$.
We remark that $\bN_{\infty,s}\big(\blalpha,\cQ\big)$ is anisotropic H\"older class of functions.

\vskip0.1cm

Let $\mathbf{E}_i,\;i=\overline{1,s}$ be the family of $s\times s$ matrices where $\mathbf{E}_i=(\mathbf{e}_1,\ldots,\mathbf{e}_i,\mathbf{0}\ldots,\mathbf{0})$ and let $\mathbf{E}_0$ is zero matrix.
Putting $K(u)=\prod_{l=1}^s\mathbf{K}(u_l),\;u_l\in\bR^s,$ we get for any $\eta\in (0,\infty)^s$ and any $z\in\bR^s$
$$
\left|\cB_{\eta,g}(z)\right|=\left|\int_{\bR^s}K(u)\left[g(z+u\eta)- g(z)\right]\rd u\right|\leq
\sum_{i=1}^s\left|\int_{\bR^s}K(u)\Big[g\big(z+\eta\mathbf{E}_{i}u\big)- g\big(z+\eta\mathbf{E}_{i-1}u\big)\Big]\rd u\right|.
$$
We note that the all components of the vectors $z+\eta\mathbf{E}_{i}u$ and $z+\eta\mathbf{E}_{i-1}u$ except $i$-th coordinate coincide.
Hence using Taylor expansion we obtain any $\eta\in (0,\infty)^s$ and  $z\in\bR^s$ in view of   (\ref{eq:inclusion-aniso-Nikolskii})
$$
\left|\int_{\bR^s}K(u)\Big[g\big(z+\eta\mathbf{E}_{i}u\big)- g\big(z+\eta\mathbf{E}_{i-1}u\big)\Big]\rd u\right|\leq \mathbf{c}Q_i\eta^{\blalpha_i}_i\int_{\bR^s}|K(u)||u|^{\blalpha_i}\rd u\leq \rk^{s}_1\mathbf{c}Q_i\eta^{\blalpha_i}_i.
$$
To get the last inequality we have taken into account (\ref{eq:ass-th:adaptation-in-supnorm}) and  used that $\mathbf{K}$ is supported on $[-1/2,1/2]$. It is worth mentioning  that  $\mathbf{c}$ as a function of $\alpha$ is bounded on any bounded domain of $(0,\infty)^s$.
Since the right hand side of the latter inequality is independent of $z$ we come to the assertion of the lemma.
\epr

\bibliographystyle{agsm}

\end{document}